\documentclass[a4paper,11pt]{amsart}
\usepackage{amssymb,amsthm,bbm,epic,eepic,graphics,amsbsy}
\usepackage{longtable}
\usepackage{a4wide}
\usepackage{amssymb,amsthm}
\usepackage{xypic}
\usepackage{array}
\usepackage{rotating}
\newcommand{\la}{\lambda}

\def\dd{\mathrm{d}}

\def\rk{\mathrm{rk}}

\def\End{\mathrm{End}}
\def\Ker{\mathrm{Ker}\,}

\def\Aut{\mathrm{Aut}}

\def\C{\ensuremath{\mathbbm{C}}}
\def\Q{\mathbbm{Q}}
\def\Z{\mathbbm{Z}}
\def\N{\mathbbm{N}}
\def\R{\mathbbm{R}}

\def\gl{\mathfrak{gl}}
\def\sl{\mathfrak{sl}}

\renewcommand{\k}{\mathbbm{k}}

\def\om{\omega}

\def\eps{\epsilon}
\def\g{\mathfrak{g}}

\def\onto{\twoheadrightarrow}
\def\into{\hookrightarrow}
\def\ii{\mathrm{i}}

\def\GL{\mathrm{GL}}

\def\bpi{\boldsymbol{\pi}}
\def\bbeta{\boldsymbol{\beta}}

\def\Id{\mathrm{Id}}

\def\Mat{\mathrm{Mat}}
\newtheorem{theo}{Theorem}[section]
\newtheorem{theointr}{Theorem}

\newtheorem{prop}[theo]{Proposition}

\newtheorem{lemma}[theo]{Lemma}
\newtheorem{remark}[theo]{Remark}

\newtheorem{conj}[theo]{Conjecture}
\newtheorem{conjintro}{Conjecture}
\newtheorem*{quest}{Question}

\def\s{\sigma}

\makeatletter
\DeclareRobustCommand\widecheck[1]{{\mathpalette\@widecheck{#1}}}
\def\@widecheck#1#2{%
   \setbox\z@\hbox{\m@th$#1#2$}%
   \setbox\tw@\hbox{\m@th$#1%
      \widehat{%
         \vrule\@width\z@\@height\ht\z@
         \vrule\@height\z@\@width\wd\z@}$}%
   \dp\tw@-\ht\z@
   \@tempdima\ht\z@ \advance\@tempdima2\ht\tw@ \divide\@tempdima\thr@@
   \setbox\tw@\hbox{%
      \raise\@tempdima\hbox{\scalebox{1}[-1]{\lower\@tempdima\box\tw@}}}%
   {\ooalign{\box\tw@ \cr \box\z@}}}
\makeatother

\def\Ad{\mathrm{Ad}}
\def\RR{\mathcal{R}}
\def\WC{W_0}

\title{{\bf Krammer representations for complex braid groups}}
\author{Ivan Marin}

\date{October 4, 2008}
\begin{document}

\maketitle

\bigskip

\begin{center}
Institut de Math\'ematiques de Jussieu\\
 Universit\'e Paris 7\\
 175 rue du
Chevaleret\\
 F-75013 Paris
\end{center}

\bigskip

\noindent {\bf Abstract.}
Let $B$ be the generalized braid group associated to some
finite complex reflection group $W$. We define a representation
of $B$ of dimension the number of reflections of the corresponding
reflection group, which generalizes the Krammer representation
of the classical braid groups, and is thus a good candidate
in view of proving the linearity of these groups. We decompose
this representation in irreducible components and compute its
Zariski closure, as well as its restriction to parabolic subgroups.
We prove that it is faithful when $W$ is a Coxeter group
of type ADE and odd dihedral types, and conjecture its
faithfulness when $W$ has a single class of reflections. If true,
this conjecture would imply various group-theoretic properties
for these groups, that we prove separately to be true for
the other groups.

\medskip


\section{Introduction}

Let $\k$ be a field of characteristic 0 and $E$ a finite
dimensional $\k$-vector space. A \emph{reflection} in $E$ is an element
$s \in \GL(V)$ such that $s^2 = 1$ and $\Ker(s-1)$ is an hyperplane
of $E$. A finite reflection group $W$ is a finite subgroup of
some $\GL(E)$ generated by a set $\RR$ of reflections. For $\k = \R$
these are finite Coxeter groups, for $\k = \C$ they are
called complex reflection groups.

To the set $\RR$ of reflections are naturally associated
an hyperplane arrangement $\mathcal{A}$, namely the set of
reflecting hyperplanes $\{ \Ker(s-1) \ | \ s \in \RR \}$,
and its complement $X = E \setminus \bigcup \mathcal{A}$.
When $\k = \C$, the generalized braid group associated to $W$ is $B = \pi_1(X/W)$. There is
a short exact sequence $1 \to P \to B \to W \to 1$,
where $P = \pi_1(X)$ is the pure braid group associated to $W$.
When $W$ is actually a finite Coxeter group, then  $B$ is an Artin group
of finite Coxeter type, and in particular for $W$ of Coxeter type $A_{n-1}$
we recover the usual braid group on $n$ strands.

This construction is also valid for (finite) pseudo-reflection groups,
where $s \in \GL(E)$ with $\Ker(s-1)$ an hyperplane
is called a pseudo-reflection if it has finite order. However,
it is known that this apparent generalization is fake, as far as group-theoretic
properties of $B$ are concerned. Indeed, using the Shephard-Todd
classification of irreducible complex pseudo-reflection groups
and case-by-case results of \cite{BMR}, we have the following fact :
for a group $B$, the following properties are equivalent
\begin{enumerate}
\item $B = \pi_1(X/W)$ for some finite pseudo-reflection group $W$. 
\item $B = \pi_1(X/W)$ for some finite reflection group $W$. 
\end{enumerate}
This was noticed for instance in \cite{BESSIS}, where indications
are given towards a possible explanation of this phenomenon (see remark 2.3
there). Recall that (finite) complex reflection groups have been classified.
According to Shephard-Todd notations, there are two
infinite series each depending on two integral parameters,
$G(2e,e,r)$ and $G(e,e,r)$, plus 15 exceptions
which include a somewhat special type called $G_{31}$ (see the appendix
for a complete list).

It is widely believed that these complex braid groups share similar properties
with Artin groups of finite Coxeter type. We refer to \cite{BESSIS} for
recent developments in this direction. For instance it is known that, if
$W$ is irreducible, then $B$ has infinite cyclic center, with the
possible (though unlikely) exception of one case ; these
groups are also torsion-free, as they have finite cohomological
dimension because $X$ is a $K(\pi,1)$ ; they have the same Garside-theoretic
flavour as the Artin groups.

Recently, an important aspect of Artin groups (of finite Coxeter
type) has been unveiled, following the work of D. Krammer and S. Bigelow
in \cite{KRAM0,KRAM,BIG},
namely that they are linear over some field of characteristic 0
(see \cite{CW,DIGNE}).
This implies that these groups are residually finite and Hopfian,
and other group-theoretic properties followed from further
investigations of this representation (see \cite{KRAMMINF}). 
The linearity result was proved by exhibiting a representation that mimics
the Krammer one for the braid groups. Characterizations
of these representations have been given, in \cite{DIGNE}
and \cite{CGW}.

The main purpose of this work is to construct a representation $R$
of $B$ that generalizes the Krammer representation. This
representation $R$ is defined over the field $K = \C((h))$ of (formal)
Laurent series and depends on an additionnal parameter $m$. It has
dimension $\# \RR$ and is
equivalent of the Krammer representation in case $W$ has simply-laced
(ADE) Coxeter type. 

\medskip

We construct this representation as the monodromy of a
$W$-equivariant integrable
1-form over $X$ of type $\om = h \sum t_s \om_s$, where $s$ runs over
$\mathcal{R}$, $\om_s$ is the logarithmic 1-form associated
to $H = \Ker(s-1) \in \mathcal{A}$, $h$ is a formal parameter
and $t_s \in \gl(V)$
where $V$ is a (finite-dimensional) complex vector space
equipped with a linear action of $W$. Such
1-forms are sometimes called generalized Knizhnik-Zamolodchikov
(KZ) systems.

We first prove (\S 2) several results in this general setting,
mostly generalizing \cite{REPTHEO}, in order
to deduce representation-theoretic properties of the monodromy representations
from the properties of the 1-form. In particular,
whenever $W_0$ is a parabolic subgroup of $W$, we may consider the integrable
1-form $\om_0 = h \sum_{s \in \mathcal{R}_0} t_s \om_s$ over the
corresponding hyperplane complement $X_0$. Generalizing results
of J. Gonz\'alez-Lorca in Coxeter type $A$, we prove for
a class of natural embeddings $\pi_1(X_0) = B_0 < B$
the following (see theorem \ref{monorestr}).

\begin{theointr} The restriction to $B_0$ of the monodromy representation
of $\om$ is isomorphic to the monodromy representation of $\om_0$.
\end{theointr}

The specific representation introduced here is associated to endomorphisms
$t_s \in \gl(V)$, where $V$ has basis $\{ v_s ; s \in \mathcal{R} \}$
and $W$ acts on $V$ by conjugation of $\mathcal{R}$.
These endomorphisms are defined by the simple formula
$$
\left\lbrace \begin{array}{lcl}
t_s.v_s & = & m v_s \\
t_s.v_u & = & v_{sus} - \alpha(s,u) v_s
\end{array} \right.
$$
where $u \neq s$ and $\alpha(s,u) = \{ y \in \mathcal{R} \ | \ ysy = u \}$.
Our central result is the following (see proposition \ref{omintegr}, theorem \ref{theorepmono}
and theorem \ref{thprops}).

\begin{theointr}
For all $m \in \C$, the 1-form $\om = h \sum t_s \om_s$
is integrable. For all but a finite set of values of $m$ :
\begin{enumerate}
\item The monodromy representation $R$ is semisimple
\item There is a natural bijection between the conjugacy classes
of reflections in $W$ and the irreducible components $R_c$ of $R$
\item For all $c$, $R_c(B)$ is Zariski-dense in the corresponding
general linear group
\item $R(P)$ is residually torsion-free nilpotent.
\end{enumerate}
\end{theointr}

Moreover, using theorem 1, we prove that this representation behaves nicely with respect
to maximal parabolic subgroups (see theorem \ref{theorepmono}).
It should be outlined that,
as far as we know, this is the first non-trivial construction
of a local system defined for arbitrary reflection groups, besides
the classical Cherednik system which gives rise to Hecke algebra
representation (see \S 4). This infinitesimal representation
is based itself on a quadratic form on the vector space
spanned by the reflections of the reflection group we start with.
The quadratic form and its basic properties are described in \S 3.
We also remark that the above theorems are `case-free', meaning
that their proof do not rely upon the Shephard-Todd classification.

\medskip

When $W$ has ADE Coxeter type, we introduced and studied this 1-form
and the corresponding monodromy representation in \cite{KRAMMINF}.
Indeed, in type ADE we have $\alpha(s,u) = 0$ if $su = us$
and $\alpha(s,u) = 1$ otherwise, thus the formula above coincides
with the one in \cite{KRAMMINF}.
We proved
there that, for generic parameters, this monodromy representation is equivalent
to the Krammer representation, which is known to be faithful. We prove
here (proposition \ref{propfaithI2}) that $R$ is also faithful when $W$ is a dihedral group of odd type.
This yields the following result.

\begin{theointr} If $W$ is a Coxeter group of type ADE or $I_2(2k+1)$,
then $R$ is faithful for generic values of $m$.
\end{theointr}

Whenever this representation $R$ is faithful, then $P$ and $B$ have the same
properties that were deduced from the faithfulness of the Krammer
representation in \cite{KRAMMINF}. It is thus a good candidate for proving the linearity of
these groups as well as generalizing properties of the
Artin groups to this more general setting.

\bigskip

In particular, this strongly suggests the following conjectures.
\begin{conjintro} If $W$ is a finite reflection group having
a single conjugacy class of reflections,
then $B$ can be embedded in $\GL_N(K)$
as a Zariski-dense subgroup, where $K$ is some field of characteristic
0 and $N$ is the number of reflections of $W$, in such a way that
$P$ is mapped into a residually torsion-free nilpotent subgroup.
In particular, $B$ is linear and $P$ is residually torsion-free nilpotent.
\end{conjintro}

\begin{conjintro} If $W$ is an irreducible finite reflection group, then $B$ can be embedded in some $\GL_N(K)$
as a Zariski-dense subgroup, where $K$ is some field of characteristic
0.\end{conjintro}

\begin{conjintro} Let $W \subset \GL_n(\C)$ be a finite pseudo-reflection group,
and $\mathcal{A}$ be the corresponding hyperplane arrangement in $\C^n$.
Then the fundamental group of the complement $\C^n \setminus \bigcup \mathcal{A}$
is residually torsion-free nilpotent.
\end{conjintro}

Conjecture 1 would be a direct consequence of theorem 2 and the faithfulness
of $R$ in the case of a single conjugacy class.
Conjectures 2 and 3 may seem less supported -- although $R$ itself
might be faithful for arbitrary $W$, and conjecture 3 is sometimes
stated for arbitrary hyperplane complements. In any case
we prove (theorem \ref{thpropsB}) that,
if conjecture 1 is true, then conjectures 2 and 3 are also true.
For this we notice that conjecture 3 is known to hold for $W$ a
Coxeter group, and prove the following
(see propositions \ref{propzarGL} and \ref{propresnil}) -- recall that a Shephard group is a complex
pseudo-reflection group which is the symmetry group of a regular complex polytope.

\begin{theointr}
Let $W \subset \GL_n(\C)$ be a finite Shephard group,
and $\mathcal{A}$ be the corresponding hyperplane arrangement in $\C^n$.
Then the fundamental group of the complement $\C^n \setminus \bigcup \mathcal{A}$
is residually torsion-free nilpotent.

\end{theointr}
\bigskip

Finally, \S 8 is an appendix devoted to the discriminant of the
quadratic form at the core of our construction. We compute
it for exceptional and Coxeter groups, plus a few other cases.

\bigskip

{\bf Acknowledgements. } This work was triggered by a discussion
with D. Krammer, who told me that one student of his, Mark Cummings, had
discovered\footnote{New Garside structures were independantly
discovered by R. Corran and M. Picantin.}
a new Garside structure on the complex braid groups
of type $G(e,e,r)$. This was indeed a stimulus for generalizing
results of \cite{KRAMMINF} to the groups $G(e,e,r)$, and then to the other
reflection groups. I also warmly thank J.-Y. H\'ee for noticing a mistake
in a previous version of this paper.

\section{Monodromy representations}

\subsection{Definitions and general facts}

Let $W \subset \GL(E)$ a finite reflection group, where $E$ is a
finite-dimensional complex vector space. We denote by $\mathcal{R}$ the
set of reflections in $W$. We let $\mathcal{A}$
denote the corresponding hyperplane arrangement in $E$, namely
the collection of $\Ker(s-1)$ for $s \in \mathcal{R}$. There
is an obvious bijection between $\mathcal{A}$ and $\mathcal{R}$,
which sends $s \in \RR$ to $\Ker(s-1)$.

Let $X = E \setminus \bigcup \mathcal{A}$ the complement of $\mathcal{A}$.
The pure braid group $P$ associated to $W$ is the fundamental
group $\pi_1(X)$, and the braid group $B$ is $\pi_1(X/W)$.
Letting $\k$ denote a field of characteristic 0, we define $\mathcal{T}$ to be the Lie algebra
over $\k$ with generators
$t_H, H \in \mathcal{A}$ and relations
$[t_Z, t_{H_0}] = 0$ for $H_0 \in \mathcal{A}$ whenever
$Z$ is a codimension 2 subspace of $E$ contained in
$H_0$, where by convention
$
t_Z = \sum_{Z \subset H} t_H$.
Note that $T = \sum_H t_H$ is obviously central in $\mathcal{T}$.

It is known by work of Kohno (see \cite{KOHNO} prop. 2.1)
that this algebra is
the holonomy Lie algebra of $X$, namely the quotient of the
free Lie algebra on $H_1(X,\k)$ by the image
of the transposed cup-product $H_2(X,\k) \to \Lambda^2 H_1(X,\k)$.
Assuming $\k = \C$ for the remainder of this section, it
follows that, given a representation $\rho : \mathcal{T} \to \gl_N(\C)$,
the closed 1-form $\om_{\rho} \in \Omega^1(X) \otimes \gl_N(\C)$
defined by
$$\om_{\rho} = \frac{1}{\ii \pi} \sum_{H \in \mathcal{A}} \rho(t_H) \om_H,
$$
where $\om_H = \mathrm{d} b_H / b_H$
for some linear form $b_H$ on $E$ with kernel $H$, satisfies
$[ \om_{\rho} \wedge \om_{\rho}] = 0$ in $H^2(X,\C)$. By a lemma of Brieskorn
(see \cite{BRIESKBOURB} lemma 5),
it follows that $\om_{\rho} \wedge \om_{\rho} = 0$, namely that
the closed 1-form $\om_{\rho}$ is integrable. Let $h$ be
a formal parameter and let us denote $A = \C[[h]]$ the ring of
formal series in $h$. Once
a base point $\underline{z} \in X$ is chosen, the equation
$\mathrm{d} F = h \om_{\rho} F$ for $F$ on $X$ with values in $A^N$
defines a monodromy representation $R : P \to \GL_N(A)$,
where $h$ is a formal parameter and $A = \C[[h]]$. We refer
to \cite{CHEN1,CHEN2} for the basic notions of formal monodromy
that are involved here.

In particular it is knwon that these monodromy representations of $P$
factorize through a universal monodromy morphism $P \to \exp \widehat{\mathcal{T}}$,
where $\widehat{\mathcal{T}}$ is the completion of $\mathcal{T}$
with respect to the graduation given by $\deg t_H = 1$. Letting
$\Mat_N(A)$ denote the set of $N \times N$ matrices with
coefficients in $A$, it follows
that $R(P) \subset \exp h \Mat_N(A)$, and in particular
$R(P) \subset 1 + h \Mat_N(A)$ ; and also that $R(P) \subset \rho(\mathsf{U}
\mathcal{T})[[h]]$ where $\mathsf{U} \mathcal{T}$
denotes the universal envelopping algebra of $\mathcal{T}$.
Letting $K = \C((h))$ be the
field of fractions of $A$, we denote by $\overline{R(P)}$
the Zariski closure of $R(P)$ in $\GL_N(K)$.

\begin{prop} \label{propintegrirred} Let $R : P \to \GL_N(K)$ be the monodromy
representation associated to $\rho : \mathcal{T} \to \gl_N(\C)$.
\begin{enumerate}
\item $\rho$ is irreducible if and only if $R$ is irreducible. In this case
$R$ is absolutely irreducible.
\item The Lie algebra of $\overline{R(P)}$ contains
$\rho(\mathcal{T}) \otimes K$.
\end{enumerate}
\end{prop}
\begin{proof}
(1) If $\rho$ is reducible and $U \subset \C^N$ is stable, then
$U \otimes_{\C} K$ is $R(P)$-stable because, for all $x \in P$,
$R(x) \in \rho(\mathsf{U} \mathcal{T}) [[h]]$. Conversely,
assume that $\rho$ is irreducible. Since $\C$ is algebraically
closed it follows by Burnside theorem that
$\rho(\mathsf{U} \mathcal{T}) =  \Mat_N(\C)$. We will prove that
$R(P)$ generates $\Mat_N(\C)$ as an associative algebra.

For $H \in \mathcal{A}$ we choose
a loop $\gamma_{H} \in \pi_1(X, \underline{z})$ around the hyperplane $H$
as follows.
Let $\Delta = \bigcup (\mathcal{A} \setminus \{ H \} )$ and
$x \in H \setminus \Delta$. Since $H$ and $\Delta$ are closed in $E$,
there exists $v \in E \setminus H$ such that $x + \la v \not\in \Delta$ 
for $|\la| \leq 1$. We define $\nu_H : [0,1] \to E$
by $\nu_H(u) = x + e^{2\ii \pi u} v$. This is a loop in $X$
with base point $x + v$. Since $X$ is connected we can choose a
path $\tau$ from $\underline{z}$ to $x+v$. Then $\gamma_H$
is defined to be the composite of $\tau$, $\nu_H$ and $\tau^{-1}$.
It is clear that $\int_{\gamma_{H_1}} \om_{H_2} = 2  \delta_{H_1,H_2}$
for $H_1,H_2 \in \mathcal{A}$.   

By Picard approximation or Chen's iterated integrals (see \cite{CHEN1,CHEN2}) we get
$$
R(\gamma_{H_0}) \in 1 + h \sum_{H \in \mathcal{A}} \rho(t_H)
\int_{\gamma_{H_0}} \om_H + h^2 \Mat_N(A) = 1 + 2h \rho(t_{H_0}) + h^2 \Mat_N(A)
$$
In particular $(R(\gamma_{H_0}) - R(1))/2 h$ belongs to
$\rho(t_{H_0}) + h \Mat_N(A)$.
Let $C$ be the $A$-subalgebra with unit of the group algebra $K P$ which
is generated
by the elements $([\gamma_{H_0}] - 1)/2  h$. We consider the algebra morphism $R : C \to \Mat_N(A)$.
It is a $A$-module morphism such whose composition with
$\Mat_N(A) \to \Mat_N(A)/h\Mat_N(A)$ is surjective. Then
by Nakayama's lemma this morphism is surjective, hence
$R(K P) = \Mat_N(K)$ and $R$ is absolutely irreducible.

We now prove (2). Let $Q_1,\dots,Q_r$ be polynomials in $N^2$
variables defining $\overline{R(P)}$ in $\GL_N(K)$. Without loss of generality,
we can assume that the $Q_i$ have coefficients in $A$. We
know that $R([\gamma_H]) = \exp h X_H$ for some $X_H \in 2 
\rho(t_H) + h \Mat_N(A)$. It follows that $\exp(m h X_H) \in R(P)$
for all $m \in \Z$. Let $u$ be a formal parameter and
consider
$Q_i(\exp(u h X_H)) = \sum_{j=0}^{\infty} Q_{i,j}^{H}(u) h^j. $
with $Q_{i,j}^H \in \C[[u]]$. It is clear that $Q_{i,j}^H \in \C[u]$.
Since $Q_{i,j}^H(m) = 0$ for all $m \in \Z$ we have $Q_{i,j}^H = 0$
for all $i,j,H$ hence $\exp(u h X_H)$ is a $K[[u]]$ point
of $\overline{R(P)}$. By Chevalley's formal exponentiation theory (see \cite{CHEVALLEY}
vol. 2 \S 8 prop. 4) it follows that $h X_H$
hence $X_H$ belongs to the Lie algebra of $\overline{R(P)}$.
On the other hand, $X_H \in \rho(\mathcal{T})[[h]] = \rho(\mathcal{T})
\otimes A$. Let $C$ be the Lie $A$-subalgebra of $\rho(\mathcal{T})
\otimes A$ generated by the elements $X_H$ and consider the
inclusion morphism. It is a $A$-module morphism whose composite
with the quotient map by $h \rho(\mathcal{T})[[h]]$ is surjective,
as $\mathcal{T}$ is generated by the $t_H$. It follows from
Nakayama's lemma that $\rho(\mathcal{T})
\otimes A$ is generated by the $X_H$, thus proving that
the Lie algebra of $\overline{R(P)}$ contains $\rho(\mathcal{T}) \otimes K$.

\end{proof}

\subsection{Monodromy of special elements}

It is known that, if $W$ is irreducible, then
$Z(P)$ is infinite cyclic and generated by
the class $\bpi$ of the loop
$u \mapsto e^{2 \ii \pi u} \underline{z}$,
except possibly if $W$ has type $G_{31}$.
This is proved in \cite{BMR}, theorem 2.24 for all
but a few cases, and \cite{BESSIS} thm. 0.5 for the remaining ones.
Since $\bpi \in P$ we have $R(\bpi) \in 1 + h \Mat_n(A)$.
If $R$ is absolutely irreducible then $R(\bpi)$
is a scalar, that is $R(\bpi) \in 1 + h A$. It follows that
$R(\bpi)$ is uniquely determined by its determinant $R(\bpi)^N$.
On the other hand, $R(\bpi) \in \exp h \rho(\mathcal{T})[[h]]$,
hence if $R(\bpi) \in \exp \left( h \rho(t) + h^2 \Mat_N(A) \right)$
for some $t \in \mathcal{T}$ then
$\det R(\bpi) = \exp h \mathrm{tr}(\rho(t))$. It follows
that $R(\bpi) = \exp( \frac{h}{N} \mathrm{tr}(\rho(t)) )$.
Moreover, it is known (see \cite{BMR}, lemma 2.4 (2)) that
$$
R(\bpi) \in 1 + 2h  \sum_{H \in \mathcal{A}} \rho(t_H) + h^2 \Mat_N(A).
$$

We thus proved the following.

\begin{prop} \label{monopi} Let $W$ be irreducible. If
$R$ is absolutely irreducible then
$$
R(\bpi) = \exp\left( \frac{2h}{N}  \sum_{H \in \mathcal{A}} \mathrm{tr}(\rho(t_H)) \right)
$$
\end{prop}

Let $H \in \mathcal{A}$, and denote $s \in \RR$
the corresponding reflection. We endow $E$ with a nondegenerate $W$-invariant
unitary form. 
Let $\Delta = \bigcup (\mathcal{A} \setminus \{ H \} )$.
Let $x \in H \setminus \Delta$, and
$v \in H^{\perp} \setminus \{ 0 \}$ such that
$x + \la v \not\in \Delta$ 
for $|\la| \leq 1$. Note that $s.v = -v$ and $s.x = x$.
We define $\nu_H : [0,1] \to E$
by $\nu_H(u) = x + e^{\ii \pi u} v$. This is a path in $X$
from $y = x + v$ to $s.y$. Let $\tau$ be a path in $X$ from
$\underline{z}$ to $y$. Then the composite of $\tau$, $\nu_H$
and $\tau^{-1}$ is a path from $\underline{z}$ to $s.\underline{z}$
in $X$, hence induces a loop in $X/W$. The class in $B$
of such a loop is called a \emph{braided reflection} around $H$.
Note that every two braided reflection around $H$ are conjugated
by an element of $P$.

We now assume that $\C^N$ is endowed with an action of $W$.
The linear action of $W$ on $E$ induces a permutation action on
$\mathcal{A}$, hence a natural action of $W$ on $\mathcal{T}$
by automorphisms of Lie algebras, where $w \in W$ maps $t_H$
to $t_{w(H)}$. Under the natural correspondance between $\mathcal{A}$
and $\RR$, the action of $W$ on $\mathcal{A}$ corresponds
to the conjugation action of $W$ on $\RR$.
The representation $\rho$
is called $W$-equivariant if $\rho(t_{w(H)}).x = w \rho(t_H) w^{-1}.x$
for all $x \in \C^N$. In this case, the representation $R$
extends to a representation $R : B \to \GL_N(K)$.

The following fact is standard.

\begin{prop} \label{monobraidref} Let $s \in \RR$, and let
$\sigma$ be a braided reflection around $H = \Ker(s-1) \in \mathcal{A}$.
Then $R(\sigma)$ is conjugated to $\rho(s) \exp( h \rho(t_H))$
in $\GL_N(K)$.
\end{prop}
\begin{proof}
Since we are only interested in the conjugacy class of $R(\sigma)$,
we can assume that $\underline{z} = x+v$ with $v \in H^{\perp}$,
$x \in H$,  $x + \la v \in X$ for all $ 0 < |\la| \leq 1$. The
differential equation $\mathrm{d }F = h\omega_{\rho} F$
along $\gamma$ has then the form
$f'(\la) = (\frac{h \rho(t_H)}{\la} + h g(\la))f(\la)$, with
$f = F \circ \gamma$ and $g$ holomorphic
in a open neighborhood of the unit disc. The conclusion follows then
from standard arguments for formal differential equations, see e.g. \cite{DIEDRAUX}
lemme 13.
\end{proof}
 
Let $\varepsilon : K \to K$ be the field automorphism sending
$f(h)$ to $f(-h)$. If $\C^N$ is endowed with a symmetric bilinear form $(\ | \ )$,
we extend it to $K^N$ and
define a skew-bilinear form $<\ | \ >$ on $K^N$
with respect to $\varepsilon$ by the formula $<x|y> = (x|\varepsilon(y))$.
Let $U_N^{\varepsilon}(K)$ denote the corresponding unitary group. We have
the following.

\begin{prop} \label{propunitaire} Assume that, for all $H \in \mathcal{A}$, $\rho(t_H)$
is selfadjoint with respect to $(\ | \ )$. Then
$R(P) \subset U_N^{\varepsilon}(K)$. Moreover, if $\om_{\rho}$
is $W$-equivariant and all $w \in W$ act orthogonally with respect to
$(\ | \ )$, then $R(B) \subset U_N^{\varepsilon}(K)$.
\end{prop}
\begin{proof}
Let $\gamma$ be a path in $X$ from $x_1$ to $x_2$. Let $F$ be solution
of $\mathrm{d} F = h \om_{\rho} F$ in a neighborhood of some
point in $\gamma([0,1])$. Let $v_1,v_2 \in K^N$ and
consider the function $g(z) = <F(z) v_1 | F(z) v_2>$ in this
neighborhood. By assumption each $\rho(t_H)$ is selfadjoint
w.r.t. $(\ | \ )$, hence $h \rho(t_H)$ is skew-symmetric with respect to
$<\ | \ >$. In particular
$$
\mathrm{d} g = <h \om_{\rho} F v_1 | Fv_2> + <Fv_1 |h \om_{\rho} F v_2 > = 0.
$$
It follows that the monodromy from $x_1$ to $x_2$ of $F$
lies in $U_N^{\varepsilon}(K)$. For $x_1 = x_2$ this means
$R(P) \subset U_N^{\varepsilon}(K)$.
For $x_2 = w.x_1$
with $w \in  W$, since $w \in O_N(\k) \subset U_N^{\varepsilon}(K)$,
this means $R(B) \subset U_N^{\varepsilon}(K)$.
\end{proof}

It is known that $B$ is generated by braided reflections (see
\cite{BMR} thm. 2.17 (1), where braided reflections are
called ``generators-of-the-monodromy'').
Moreover, the exact sequence $1 \to P \to B \to W \to 1$
induces an exact sequence between the centers
$1 \to Z(P) \to Z(B) \to Z(W) \to 1$, and $Z(B)$ is generated
by the class $\bbeta$ of the loop $u \mapsto \underline{z}\exp(2 \ii \pi u/\# Z(W))$ (see \cite{BMR} thm. 2.24 and \cite{BESSIS} thm. 0.5),
except possibly if $W$ has type $G_{31}$.
Assume that $Z(W)$ acts trivially on $\C^N$. For all $w \in W$,
if $\tilde{w} \in B$ is in the preimage of $w$ then
$R(\tilde{w})$ acts like $w$ on $(A/hA)^N$. It follows that,
if $Z(W)$ acts trivially on $\C^N$ and $R$ is absolutely irreducible,
we have $R(\bbeta) \in 1 + h A$. Since $\bbeta^{\# Z(W)} = \bpi$,
and elements in $1 + h A$ have unique roots of any order,
we proved the following.

\begin{prop} \label{monobeta} Assume that $W$ is irreducible and that $\C^N$
is endowed with a linear action of $W$ such that $Z(W)$ acts trivially.
If $\rho$ is $W$-equivariant and the restriction of $R$ to $P$
is absolutely irreducible,
then
$$
R(\bbeta) = \exp\left( \frac{2h}{N \#Z(W)}  \sum_{H \in \mathcal{A}} \mathrm{tr}(\rho(t_H)) \right)
$$
\end{prop}

\begin{remark} \label{remcentre}
\end{remark}
Let $\rho_0 : \mathcal{T} \to \gl_1(\C)$ be defined by
$\rho_0(t_s) = \Id$. This is $W$-equivariant with respect to
the trivial action of $W$ on $\C$. Let $R_0 : B \to \GL_1(K)$ be
the corresponding monodromy representation. By proposition \ref{monobraidref},
the image of any braided reflection is $q = \exp( h)$. It is known
by \cite{BMR} that $B$ is generated by braided reflections, hence $R_0$
factors through a morphism $\varphi : B \to \Z$ which sends braided reflections
to $1$. By the proposition we have
$\varphi(\bbeta) = 2  \# \mathcal{R} /\# Z(W)$.

\bigskip

Instead of deducing proposition \ref{monobeta} from proposition \ref{monopi},
we could also have deduce it from the following one.

\begin{prop} \label{monoregu} Let $w \in W$, $\underline{z} \in X$ such that $w.\underline{z} = e^{\ii \theta} \underline{z}$
for some $\theta \in ]0,2 \pi]$. The path $\gamma : u \mapsto e^{\ii u \theta}$
takes its values in $X$ and its image $[\gamma] \in \pi_1(X/W,\underline{z})$ satisfies
$$
R([\gamma]) = w \exp\left( h \frac{\theta}{\pi} \sum_{H \in \mathcal{A}} \rho(t_H)\right).
$$
\end{prop}
\begin{proof}
For all $H \in \mathcal{A}$ we have $b_H(e^{\ii u \theta} \underline{z}) = b_H(\underline{z}) \neq 0$
hence $\gamma([0,1]) \subset X$. 
From $\om_H = \dd b_H/b_H$
we get $\gamma^* \om_H = (\ii \theta) \dd u = (\theta/\pi) \dd u$.
Since $\om_{\rho} = (1/\ii \pi) \sum \rho(t_H) \om_H$ we have
$\gamma^* \om = (\theta/\pi) (\sum \rho(t_H))
\dd u$ hence $f(u) = \exp(uh (\theta/\pi) \sum \rho(t_H))$ satisfies
$\dd f = h (\gamma^* \om_{\rho})f$
and $f(0) = \mathrm{Id}$. From this we get $R([\gamma]) = w f(1)$
and the conclusion.
\end{proof}

\subsection{Restriction to parabolic subgroups}

Let $I \subset E$ be a proper subspace. The parabolic subgroup $W_0$
of $W$ associated to $I$ is the subgroup of $W$ of the
elements which stabilize $I$ pointwise. By a result of Steinberg (see
\cite{STEINBERG}, thm. 1.5),
it is generated by the set $\RR_0$ of reflections of $W$ whose reflecting
hyperplanes contain $I$, and can be considered as a reflection
group acting either on $E$ or $I^{\perp}$. We let $\mathcal{A}_0
= \{ H \in
\mathcal{A} \ | \ H \supset I \}$, and $X_0 = E \setminus \bigcup \mathcal{A}_0$.
Notice that any such $W_0$ is the starting point of a (finite) chain
$W_0 < W_1 < \dots < W$ of reflection groups where each $W_i$ is a \emph{maximal}
parabolic subgroup of $W_{i+1}$. For maximal parabolic subgroups, $I$ is
a complex line.

Let $\mathcal{T}_0$ be the holonomy Lie algebra associated to such a
parabolic subgroup $W_0 \subset
\GL(E)$. The following lemma shows that we can identify $\mathcal{T}_0$
with a sub-Lie-algebra of $\mathcal{T}$.

\begin{lemma} \label{embedLie} The inclusion $\mathcal{R}_0 \subset \mathcal{R}$
induces a $W$-equivariant Lie algebra embedding $\mathcal{T}_0
\to \mathcal{T}$ which sends $t_H$ to $t_H$ for $H \in \mathcal{A}_0$.
\end{lemma}
\begin{proof}
Without loss of generality we can assume that $W_0$ is maximal among parabolic subgroups.
Recall that $W_0 = \{ w \in W \ | \ \forall x \in I \ \ w.x = x \}$.
Here we denote $t_H^0$ the generators of $\mathcal{T}_0$ for
$H \in \mathcal{A}_0$. If $Z = H_1 \cap H_2$ is a codimension 2 subspace
of $E$ with $H_1,H_2 \in \mathcal{A}_0$ then $I \subset Z$, and
$$
t_Z^0 = \sum_{\stackrel{H \in \mathcal{A}_0}{H \supset Z}} t_H^0 
= \sum_{\stackrel{H \in \mathcal{A}}{H \supset Z}} t_H.
$$
This proves that $t_H^0 \mapsto t_H$ can be extended (uniquely)
to a Lie algebra morphism $j : \mathcal{T}_0 \to \mathcal{T}$,
which is obviously $W$-equivariant. In order to
prove that $j$ is injective, we define a Lie algebra morphism $q : \mathcal{T} \to
\mathcal{T}_0$ such that $q \circ j$ is the identity of $\mathcal{T}_0$.

We define $q$ by sending $t_H$ to $t_H^0$ is $H \in \mathcal{A}_0$,
and to $0$ otherwise. In order to prove that it is
well-defined, we only need to check that, if $Z$ is a codimension 2
subspace of $E$ not containing $I$ and $H \supset Z$ belongs to
$\mathcal{A} \setminus \mathcal{A}_0$, then $[q(t_Z),q(t_H)] = 0$,
where $q(t_Z) = \sum_{H \supset Z} q(t_H)$,
all other cases being trivial. But then $H$ is the only hyperplane
in $\mathcal{A}_0$ containing $Z$, because if $H' \in \mathcal{A}_0$
were another hyperplane containing $Z$ then $Z = H \cap H'$
hence $Z \supset I$. It follows that $q(t_Z) = q(t_H) = t^0_{H}$
hence $[q(t_Z),q(t_H)] = 0$. Then $q$ is well-defined,
$q \circ j(t_H^0) = t_H^0$ for all $H \in \mathcal{A}_0$ hence
$q \circ j$ is the identity of $\mathcal{T}_0$.

\end{proof}

Let $B_0$ and $P_0$ denote the braid group and pure braid group
associated to $W_0$. It is clear that we have natural identifications
$P_0 = \pi_1(X_0)$ and $B_0 = \pi_1(X_0/W_0)$. Following
\cite{BMR} \S 2D, we define embeddings of $P_0$ and $B_0$ in $P$ and $B$,
respectively. The images of such embeddings are called
parabolic subgroups of $P$ and $B$, respectively.

We endow $E$ with a $W$-invariant unitary form and denote $\| \ \|$ the associated
norm. Let $x_1 \in I$ such that
$x_1 \not\in H$ for all $H \in \mathcal{A} \setminus \mathcal{A}_0$.
There exists
$\eps > 0$ such that, for all $x \in E$ with $\| x -x_1 \| \leq \eps$,
$x \not\in H$ for all $H \in \mathcal{A} \setminus \mathcal{A}_0$.
We let $x_2$ fulfilling $\| x_2 -x_1 \| < \eps$ and such
 that $x_2 \not\in H$
for all $H \in \mathcal{A}$.
Let $\Omega = \{ x \in E \ | \ \| x -x_1 \| \leq \eps \}$.
It is easily checked that
$\pi_1(X \cap \Omega,x_2) \to \pi_1(X_0,x_2)$ is an isomorphism,
hence the obvious inclusion $\pi_1(X \cap \Omega,x_2) \to \pi_1(X,x_2)$
defines an embedding $P_0 \to P$. Since $\Omega$ is setwise stabilized
by $W_0$, this embedding extends to an embedding $B_0 \to B$.
It is proved in \cite{BMR} \S 2D that such embeddings are well-defined up to
$P$-conjugation.

In the remainder of this section we study the restriction
of the monodromy representation $R$ to such a parabolic subgroup.
For this, we can obviously assume that $W_0$ is maximal, that is $\dim_{\C} I = 1$.
The analytic part of the argument is essentially the same
as for Knizhnik-Zamolodchikov systems and classical braid groups. Since
this case has already been dealt with in full details
by J. Gonz\'alez-Lorca in \cite{JORGE} (part 2), we
allow ourselves to be somewhat sketchy on the analytic justifications,
and focus instead on the topological and algebraic changes that are needed
for the general case.

\subsubsection{Contruction of tubes}

Let $e \in I \setminus \{ 0 \}$. For all $H \in \mathcal{A} \setminus \mathcal{A}_0$,
we have $b_H(e) \neq 0$. Since the linear forms $b_H$
are defined only up to some nonzero scalar, we can assume
$b_H(e) = 1$ for all $H \in \mathcal{A} \setminus \mathcal{A}_0$.
We can choose $x_1 \in I$ and $\Omega$ as above such that, for all $x \in \Omega$
and all $H \in \mathcal{A} \setminus \mathcal{A}_0$, we
have $b_H(x) \not\in \R^-$. It follows that,
for all $\la \in \R^+$ and all $x \in (\Omega \cap X)\setminus I$,
we have $x + \la e \in X$.

\begin{figure}
\resizebox{10cm}{!}{\includegraphics{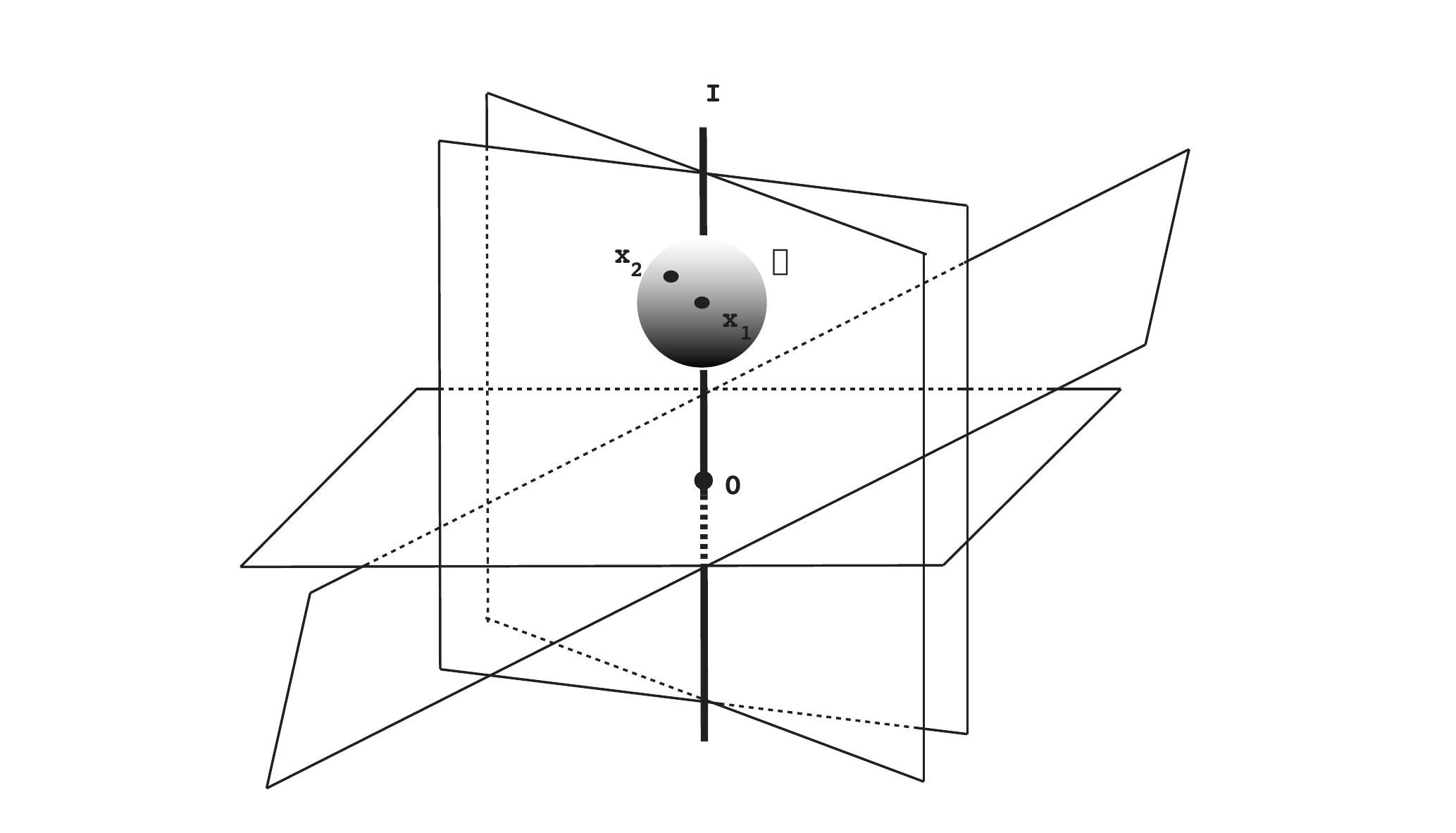}}
\caption{The ball $\Omega$}
\end{figure}

\begin{figure}
\resizebox{16cm}{!}{\includegraphics{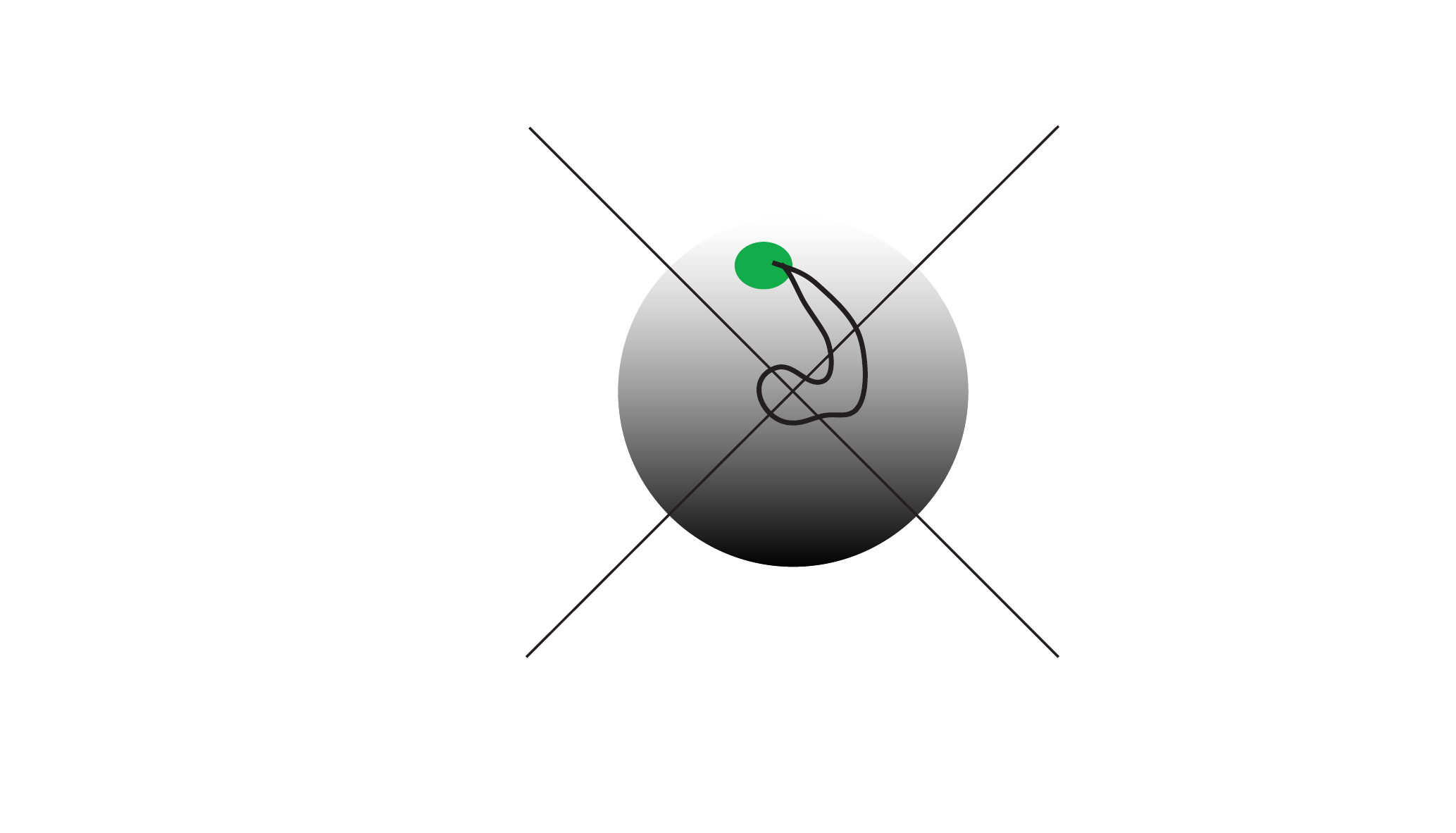}\includegraphics{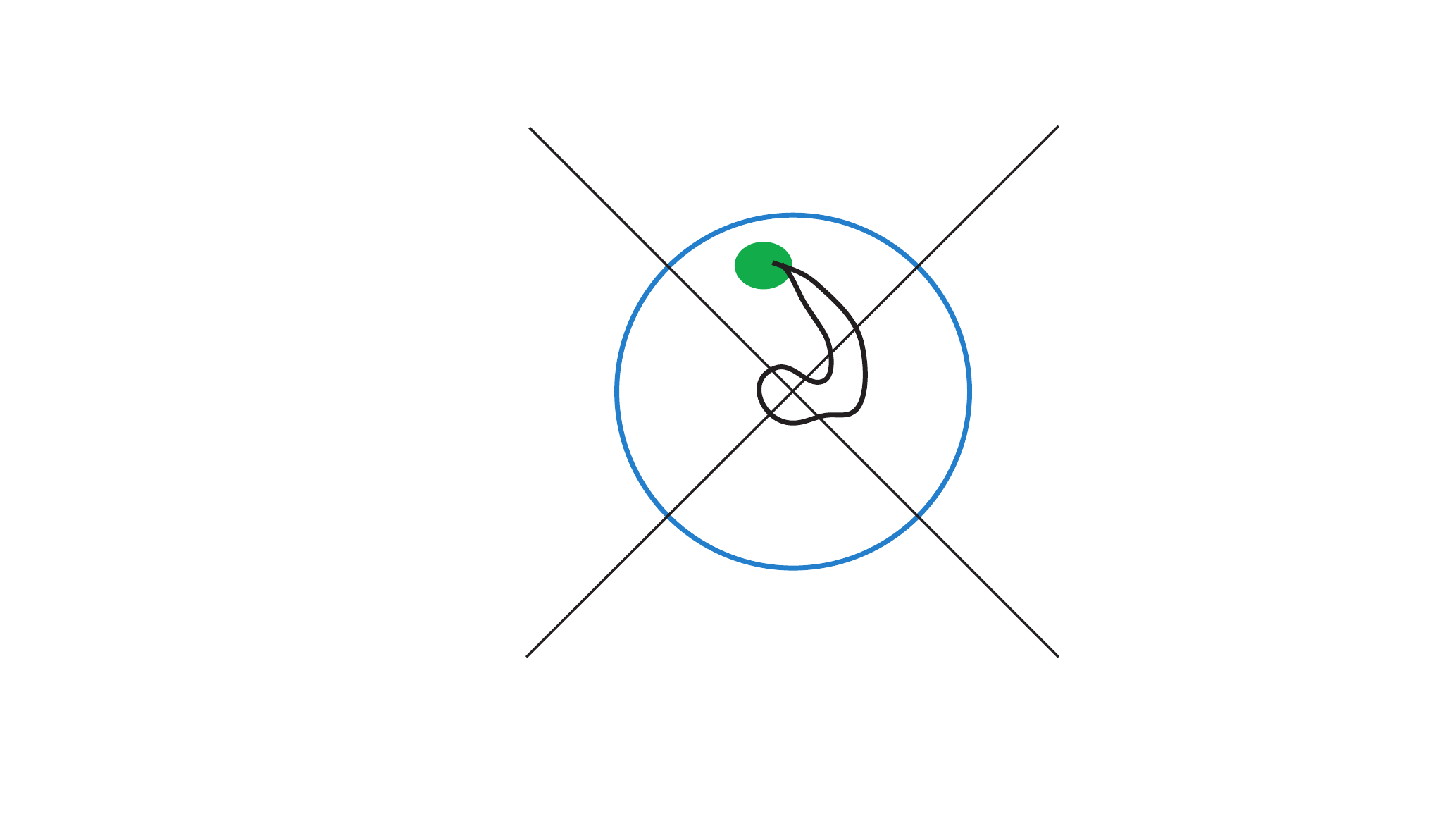}}
\caption{Monodromy in $\Omega$ or $M$}
\end{figure}

In particular, there exists a contractible open neighborhood $U$ of $x_2$ in $E$
such that $U + \R^+ e =U$. Indeed, let
$\mathcal{B}$ be an open ball centered at $x_2$ inside $(\Omega\cap X) \setminus I$.
Then $U = \mathcal{B} + \R^+ e$ is convex, open, included in
$(\Omega\cap X) \setminus I$, and satisfies $U=U+\R^+ e$.
Similarly, the open set $M = \Omega^{\circ} \cap X + \R^+ e$ contains $U$
and is homotopically equivalent to $\Omega \cap X$.
The restriction of $R$ to $P_0$ or $B_0$ can thus been considered
as the monodromy of $\mathrm{d} F = h \om_{\rho} F$ in $M$,
where $F$ has values in $\Mat_N(A)$ and $F(x_2) = \mathrm{Id}$.

\subsubsection{Restriction of the differential form}

Let $\om^0_{\rho} = \sum_{H \in \mathcal{A}_0} \rho(t_H) \om_H$ and $n+1 = \dim E$.
We choose coordinates $u,w_1,\dots,w_n$ in $E$ such that $u$ is the
coordinate along
$I$ corresponding to $e$ and $w_1,\dots ,w_n$ are coordinates in $I^{\perp}$.
Then $h\omega_{\rho}$ can be written $\theta_0 \mathrm{d}u + \theta_1 \mathrm{d} w_1
+ \dots + \theta_n \mathrm{d} w_n$. 
The integrability condition
implies in particular
$$
\frac{\partial \theta_i}{\partial u} -
\frac{\partial \theta_0}{\partial w_i} + [\theta_i,\theta_0] = 0.
$$
for $1 \leq i \leq n$.
Similarly, we can write
$\om_{\rho}^0 = \widehat{\theta_1} \mathrm{d} w_1
+ \dots + \widehat{\theta_n} \mathrm{d} w_n$, since
$H \supset I$ implies $b_H(I) = 0$, hence $b_H$
is a linear combination of $w_1,\dots,w_n$
if $H \in \mathcal{A} \setminus \mathcal{A}_0$.

On the other hand $H \not\supset I$ implies $b_H(I) \neq \{ 0 \}$.
Let $v \in I \setminus \{ 0 \}$. Then, for all $x \in E$,
$|b_H(x+\la v) | \to \infty$ when $\la \to \infty$,
and in particular, for $w_1,\dots,w_n$ fixed and $u = \la$,
$\theta_i - \widehat{\theta}_i \to 0$ when $\la \to \infty$.
Furthermore, $\theta_i - \widehat{\theta}_i = \Delta_i/u$, with
$\Delta_i$ a
bounded analytic function on $M$.

The equation $\mathrm{d} F = h\om_{\rho} F$ can thus be written as
the system of partial differential equations
$$
\frac{\partial F}{\partial u} = \theta_0 F, 
\frac{\partial F}{\partial w_i} = \theta_i F, 1 \leq i \leq n.
$$
Let $F$ be the solution in $U$ of $\mathrm{d} F = h\om_{\rho} F$
such that $F(x_2) = \mathrm{Id}$.
Assume that we have a solution $f$ in $U$ of
$\frac{\partial f}{\partial u} = \theta_0 f$ such that $f \equiv 1$ modulo $h$.
Since $F$ also satisfies
this equation and $f$ is invertible, there
exists a well-defined function $G$ with values in $\Mat_N(A)$,
independant of $u$, such that $F = fG$.
We will find a condition on $f$ ensuring that $G$ satisfies
$\mathrm{d} G = h \om_{\rho}^0 G$. It is readily checked that, under $\mathrm{d}
F = h\om_{\rho} F$,
$\partial G/\partial w_i = \widehat{\theta}_i G$
is equivalent to
$$\widehat{\theta}_i = f^{-1}\left(\theta_i f - \frac{\partial f}{\partial
w_i} \right).
$$
Let $g_i = f^{-1}(\theta_i f - \partial f/\partial
w_i )$. We find
$$
\frac{\partial g_i}{\partial u} = f^{-1} \left( 
\frac{\partial \theta_i}{\partial u} -
\frac{\partial \theta_0}{\partial w_i} + [\theta_i,\theta_0] 
\right) f = 0
$$
by the integrability assumption on $\om$. It follows that
$g_i$ is independant of $u$, hence we only need to
find $f$ such that, for fixed coordinates $w_1,\dots,w_n$,
and $u = u_0 + \la$ for $\la \in \R^+$,
$$
\widehat{\theta}_i = \lim_{\la \to +\infty }
f^{-1}\left(\theta_i f - \frac{\partial f}{\partial
w_i} \right)
$$

\subsubsection{Fuchsian differential equation along $I$}

Let $x \in M$ and consider the function $f_x : \la \mapsto f(x+\la e)$. It is
defined for $\la \in \R^+$, and the partial differential equation satisfied
by $f$ implies that $f_x$ satisfies a differential
fuchsian equation
$$
f_x'(\la) = h \left( \sum_{H \not\supset I } \frac{\rho(t_H)}{b_H(x+\la e)} 
\right) f_x(\la) = h
\left( \sum_{H \not\supset I } \frac{\rho(t_H)}{\la + b_H(x)}
\right) f_x(\la)
$$
Let $X = \sum_{H \not\supset I } t_H$. Recall that $b_H(x) \not\in \R^-$.
It is classical (see e.g. \cite{JORGE}, or \cite{DIEDRAUX} appendix 1)
that there exists a unique $f_x$ satisfying this equation such that
$f_x \sim \la^{h \rho(X)}$ when $\la \to + \infty$ with $\la \in \R$,
with $f_x$ real analytic in $]0,+\infty[$, and $f_x \sim \la^{h \rho(X)}$
meaning that there exists $g_x$ such that $f_x = \la^{h \rho(X)}(1+\frac{1}{\la} g_x)$
with $g_x$ bounded when $\la \to +\infty$. Moreover, we have $f_x \equiv 
\mathrm{Id}$
modulo $h$.

We recall the main steps of the argument.
Let $j(\la) = f_x(1/\la)$. It is defined on $]0,+\infty[$ and
satisfies
$$
j'(\la) = h \left( \frac{-\rho(X)}{\la} + \sum_{H \not\supset I}
\frac{ \rho(t_H)}{\la + \frac{1}{b_H(x)}} \right) j(\la).
$$
Then standard arguments (see e.g. \cite{DIEDRAUX} lemme 13) show
that there exists a unique solution $j$ such that
$j(\la) \sim \la^{-h\rho(X)}$ when $\la \to 0$, meaning 
$j(\la) = \la^{-h\rho(X)}(\Id + k(\la))$, for some $k$ analytic on $[0,+\infty[$
with $k(0)
= 0$. This implies $f(\la) = \la^{h \rho(X)}(1+k(1/\la))$.
Since $k(0) = 0$ and $k$ is analytic we get $k(1/\la) = g(\la)/\la$
for some analytic $g$ on $]0,+\infty[$ which remains bounded when $\la \to +\infty$.
Moreover, we have $f'_x \equiv 0$ modulo $h$, hence
$f_x$ modulo $h$ is constant. Since $k(0) = 0$ we get $j(\la) \equiv \mathrm{Id}$
modulo $h$ hence $f_x(\la) \equiv \mathrm{Id}$ modulo $h$.
We leave to the reader the verification that $f_x$ and $g_x$ vary analytically
in $x \in M$, either by checking the constuction of $j_x$ in
the proof of \cite{DIEDRAUX} lemme 13, or by using the
more explicit description in \cite{JORGE} of $j$ in terms of
Lappo-Danilevskii polylogarithms, which vary analytically in their parameters
$b_H(x)$, for $H \in \mathcal{A} \setminus \mathcal{A}_0$.

\subsubsection{Conclusion}

We now prove that $X$ commutes to all $t_{H_0}$ for $H _0\in \mathcal{A}_0$.
We remark that
$$
X = \sum_{H \in \mathcal{A}} t_H - \sum_{H \in \mathcal{A}_0} t_H
= T - T_0
$$
where $T$ and $T_0$ are central elements of $\mathcal{T}$ and
$\mathcal{T}_0$, respectively (recall that we identified $\mathcal{T}_0$
with a sub-Lie-algebra of $\mathcal{T}$ by lemma \ref{embedLie}).
Then $[t_{H_0},T] = 0$ and $[t_{H_0},T_0] = 0$ hence $[t_X,t_{H_0}] = 0$

We are now ready to prove that $f$ satisfies what needed.
Let $x_0 \in U$, with coordinates
$w_1,\dots,w_n,u_0$, and consider $x$ in $U$ with
coordinates $w_1,\dots,w_n,u = u_0 + \la$ for $\la \in \R^+$.
On the one hand, since $f = \la^{h \rho(X)}(1 + \frac{1}{\la} g)$,
we have
$$
f^{-1} \frac{\partial f}{\partial w_i} = \left(1+ \frac{1}{\la} g
\right)^{-1} \frac{1}{\la} \frac{\partial g}{\partial w_i} \to 0
$$
when $\la \to +\infty$. On the other hand, 
we have $\la^{-h\rho(X)} \theta_i \la^{h\rho(X)} =
\la^{-h\rho(X)} (\theta_i -\widehat{\theta}_i) \la^{h\rho(X)}
+ \la^{-h\rho(X)} \widehat{\theta_i} \la^{h\rho(X)}$.
We know that $\theta_i -\widehat{\theta}_i \to 0$, and more
precisely that $\theta_i -\widehat{\theta}_i = \frac{1}{\la}\widetilde{\Delta}_i(\la)$
for some bounded function $\widetilde{\Delta}_i$.
It follows that $\la^{-h\rho(X)} (\theta_i -\widehat{\theta}_i) \la^{h\rho(X)} \to 0$
when $\la \to + \infty$. Since $\rho(X)$ commutes with $\widehat{\theta}_i$,
we have $\la^{-h\rho(X)} \widehat{\theta_i} \la^{h\rho(X)}=
\widehat{\theta_i}$ for
all $\la$. It follows that
$$
f^{-1} \theta_i f = \left( 1 + \frac{1}{\la} g \right)^{-1}
\la^{-h\rho(X)} \theta_i \la^{h\rho(X)}
\left( 1 + \frac{1}{\la} g \right) \to \widehat{\theta}_i
$$
when $\la \to + \infty$.

We thus proved that $F$ can be written on $U$ as $F = f G$,
with $G$ a function independant of $u$
satisfying $\mathrm{d} G = h\om_{\rho}^0 G$. Since $f$
is analytic in $M$, it follows that $F$ and $G$ have the same
monodromy with respect to a loop in $M$, which is a deformation
retract of $X \cap \Omega$. Moreover, if $\om_{\rho}$
is $W$-equivariant, for each $w \in W_0$
the function $x \mapsto f(w.x)$ satisfies
the same differential equation, as well as the same
asymptotic conditions, as $f$. It follows that $f(w.x) = f(x)$
for all $x \in M$, hence $F$ and $G$ have the same
monodromy with respect to every path with endpoints in $W_0.x_2$.

Since, for any parabolic subgroup $W_0$ of $W$, the embeddings of $P_0$
and $B_0$ in $P$ and $B$ are well-defined up to conjugacy, and since there
exists a chain $W_0 < W_1 < \dots < W$ of parabolic subgroups with $W_i$
maximal in $W_{i+1}$, we thus proved the following.

\begin{theo} \label{monorestr} Let $W_0$ be a parabolic subgroup of $W$.
Let $\rho_0$ be the restriction of $\rho$ to $\mathcal{T}_0$, and $R_0$
be the monodromy of $\rho_0$. Then the restriction of $R$
to $P_0$ is isomorphic to $R_0$. If $\om_{\rho}$ is $W$-equivariant,
then the restriction of $R$ to $B_0$ is isomorphic to $R_0$.
\end{theo}

\begin{remark}
\end{remark}
This result is stated for a given infinitesimal representation because
this is the statement we need here. Another way to state it,
which may be useful in other contexts, is to say that the
following diagram commutes
$$
\xymatrix{
B  \ar[r] & W \ltimes \exp \widehat{\mathcal{T}} \\
B_0 \ar[u] \ar[r] & W_0 \ltimes \exp \widehat{\mathcal{T}_0} \ar[u] \\
}
$$
where the horizontal maps are the universal monodromy
morphisms associated to $x_2$ and the embedding $\mathcal{T}_0 \to
\mathcal{T}$ is defined by lemma \ref{embedLie}. Indeed, the
proof given here apply verbatim to the differential equation
with values in $\widehat{\mathsf{U} \mathcal{T}}$ given by
$\mathrm{d} F = \om F$.

\begin{remark}
\end{remark}
The definition of the holonomy Lie algebra $\mathcal{T}$
and the monodromy construction can be made for $W$ a \emph{pseudo-}reflection
group. Most of the results and
proofs of this section remain valid in this setting, including propositions \ref{propintegrirred}, \ref{monopi}, \ref{propunitaire}, \ref{monobeta}, \ref{monoregu},
lemma \ref{embedLie} and theorem \ref{monorestr}. Only the definition of reflections and braided reflections
need additional care in this more general setting.

\section{The quadratic permutation module $V$}
Let $\k$ be a field of characteristic 0 and let $W$
be a finite complex reflection group. We let $\RR$
denote its set of reflections, and $V$ a $\k$-vector space of
dimension $\# \mathcal{R}$ with basis $(v_s)_{s \in \RR}$.
It is a permutation $W$-module under the action
$w.v_s = v_{wsw^{-1}}$. The goal of this section
is to endow it with a family of $W$-invariant quadratic forms.

For $s \in \RR$ let $H_s = \Ker(s-1)$ denote the corresponding reflection hyperplane in $E$.
Every codimension 2 subspace $Z$ in $E$ can be written as $H_s \cap H_u$
for some $s,u \in \RR$ and $u \neq s$. Conversely, if $s \neq u$ then
$H_s \cap H_u$ has codimension 2. We denote $\RR_{s,u}$ for $s\neq u$
the set of reflections of $W$ which pointwise stabilize $H_s \cap H_u$.
We have $\RR_{s,u} = \RR_{u,s}$, and $w \RR_{s,u} w^{-1} = 
\RR_{wsw^{-1},wuw^{-1}}$.

We begin with a lemma.

\begin{lemma} Let $u,y \in \RR$. If $yuy \neq u$ then $y \in \RR_{u,yuy}$.
\end{lemma}
\begin{proof}
Let $Z = H_u \cap H_{yuy}$. We first show that $Z$ is setwise
stabilized by $y$. For $z \in Z$ we have $yuy.z = z$ hence
$u(y.z) = y.z$ hence $y.z \in H_u$ ; likewise $(yuy)(y.z) = yu.z = y.z$
since $z \in H_u$ hence $y.z \in H_{yuy}$ (in other words
$y.H_u = H_{yuy}$ and $y.H_{yuy} = H_u$). In particular $Z$ is
setwise stabilized by $y$, as well as its orthogonal $Z^{\perp}$
with respect to some $W$-invariant unitary form on $E$.
Since $y \in \RR$, either $Z$ or $Z^{\perp}$ is pointwise stabilized
by $y$. If $Z^{\perp}$ were pointwise stabilized by $y$, writing
$y$ and $u$ as blockmatrices on $Z \oplus Z^{\perp}$
we would have $y = \begin{pmatrix} * & 0 \\ 0 & 1 \end{pmatrix}$
and $u = \begin{pmatrix} 1 & 0 \\ 0 & * \end{pmatrix}$
hence $yuy = u$, which has been excluded. It follows that $Z$ is pointwise stabilized
by $y$, that is $Z \subset H_y$ and $y \in \RR_{u,yuy}$.
\end{proof}

For $s \neq u$ with $s,u \in \RR$ we let $\alpha(s,u) = \# \{y \in \RR_{s,u} \ | \
yuy = s\}$. By the above lemma we have
$$
\alpha(s,u) = \# \{ y \in \RR \ | \ yuy = s \}.
$$
Moreover
$$
\alpha(u,s) = \# \{ y \in \RR \ | \ yuy = s \} = 
\# \{ y \in \RR \ | \ u = ys y\} = \alpha(s,u)
$$
and, for all $w \in W$,
$$
\begin{array}{lcl}
\alpha(wsw^{-1},wuw^{-1}) & = &
\# \{ y \in \RR \ | \ ywuw^{-1}y = ws w^{-1}\}\\
&=&\# \{ y \in w\RR w^{-1} \ | \ (w^{-1}yw)u(w^{-1}yw) = s \} \\
&=&\# \{ y \in \RR \ | \ yuy = s \}\\
& =& \alpha(s,u).
\end{array}
$$
Let $c \in \RR/W$. We define a nonoriented graph $\mathcal{G}_c$ on the set $c$
with an edge between $u,s \in \RR$ iff $u \neq s$ and $\alpha(s,u) > 0$.

\begin{lemma} \label{lemGconn}
The graph $\mathcal{G}_c$ is connected.
\end{lemma}
\begin{proof}
Let $s,u \in c$ such that $s \neq u$. Since $c$ is a conjugacy
class there exists $w \in W$ such that $s = wuw^{-1}$.
Since $W$ is generated by $\RR$ there exist
$y_1,\dots, y_r \in \RR$ with $r \geq 1$ such that
$w = y_1\dots y_r$, hence $s = y_r \dots y_1 u y_1 \dots y_r$.
We can assume that $r$ is minimal with respect to this property.
This implies that $y_1uy_1 \neq u$, $y_2 y_1 u y_1 y_2 \neq y_1 u y_1$,
etc. It follows that $\alpha(u,y_1uy_1) > 0$, $\alpha(y_1uy_1,
y_2y_1uy_1y_2) > 0$, \dots, $\alpha(y_{r-1}\dots y_1 u y_1\dots y_{r-1},s)
>0$, which means that $u$ is connected to $s$ by a path in $\mathcal{G}_c$,
hence $\mathcal{G}_c$ is connected. 

\end{proof}

To $c \in \RR/W$ we associate a real matrix $A_c$ of order $\# c$ as follows.
Choose a total ordering on $c$, and let the $(s,u)$ entry
of $A_c$ be $\alpha(s,u)$ if $s\neq u$, and $1$ otherwise. This
matrix depends on the choice of the ordering, but its eigenvalues
obviously do not. Let $V_c$ be the subspace of $V$
spanned by the $v_s, s \in c$. Let $(\ ,\ )$ be the $W$-invariant quadratic
form on $V_c$ defined by $(v_s,v_u) = \delta_{s,u}$, and $N(c) \in \N$
be defined by
$$
N(c) = \frac{(A_c v_c ,v_c)}{(v_c,v_c)}
$$
where $v_c = \sum_{s \in c} v_s$. This quantity does not depend on
the chosen ordering either. Using the $W$-invariance of $\alpha$,
a straightforward computation shows that, for all $s \in c$,
$$
N(c) = 1 + \sum_{u \in c \setminus \{ s \}} \alpha(s,u) \in \N.
$$

\begin{lemma} The greatest eigenvalue of $A_c$ is $N(c)$,
and occur with multiplicity 1.
\end{lemma}
\begin{proof} For simplicity we note $A = A_c$. Let $\beta \in \R$ be the greatest eigenvalue
of $A$, and $F = \Ker(A - \beta)$. Since $A$ is symmetric, we know that $\forall x \in V_c\ \ 
(Ax,x) \leq \beta(x,x)$ and $(Ax,x) = \beta(x,x) \Leftrightarrow x \in F$.

Let $x \in F \setminus \{ 0 \}$. We write $x = \sum_{s \in c} \la_s v_s$
with $\la_s \in \R$. By convention here we let $\alpha(s,s) = 1$.
We have $(Ax,x) = \sum_{(s,u)\in c^2} \la_s \la_u \alpha(s,u)$.
Choose $s_1,s_2 \in c$ with $s_1 \neq s_2$. If $\alpha(s_1,s_2) \neq 0$, then $\la_{s_1}$ and $\la_{s_2}$
have the same sign. Indeed, since $\alpha(s_1,s_2) > 0$,
$\la_{s_1} \la_{s_2} <0$ would imply $\la_{s_1} \la_{s_2}
\alpha(s_1,s_2) < |\la_{s_1}| |\la_{s_2}| \alpha(s_1,s_2)$.
Denoting $|x| = \sum |\la_s| v_s$, we would have
$(Ax,x) < (A|x|,|x|)$ and $(x,x) = (|x|,|x|)\neq 0$ hence $\beta
< (Ax,x)/(x,x)$, a contradiction. Let now $y \in \RR$
such that $s_2 = ys_1 y$. Since the quadratic form
$z \mapsto (Az,z)$ is $W$-invariant, the vector space
$F$ is setwise stabilized by $W$. If follows that $y.x -x \in F$.
But $y.x-x = (\la_{s_2}- \la_{s_1})v_{s_1} + 
(\la_{s_1}- \la_{s_2})v_{s_2}+ \dots$, hence $y.x-x \in F$
implies that $\la_{s_1}- \la_{s_2}$ and $\la_{s_2}- \la_{s_1}$
have the same sign, meaning $\la_{s_2}= \la_{s_1}$. It follows that
$s \mapsto \la_s$ is constant on each connected component of
$\mathcal{G}_c$. Since $\mathcal{G}_c$ is connected, this
means that $F$ is generated by $v_c$, whence the conclusion.
\end{proof}

We define a symmetric bilinear form on $V$ depending on $m \in \k$
by $(v_s|v_u) = \alpha(s,u)$ if $s \neq u$ and $(v_s|v_s) = (1-m)$.
Since $\alpha(s,u) = 0$ if $s,u$ belong to different classes
in $\RR/W$, the direct sum decomposition
$$
V = \bigoplus_{c \in \RR/W} V_c
$$
is orthogonal w.r.t. the form $(\ | \ )$. In particular,
this form restricts to symmetric bilinear forms $(\ | \ )_c$
on each $V_c, c \in \RR/W$. It is clear that the matrix $A_c - m$
is the matrix of $(\ | \ )_c$ with respect to a basis of $V_c$ made
out of vectors $v_s$ for $s \in c$. By abuse of terminology, we call $\det(A_c - m)$
the \emph{discriminant} of $(\ | \ )_c$. For irreducible Coxeter groups of type $ADE$,
this form gives the (opposite of) the one introduced
in \cite{KRAMMINF}.

\begin{prop} \label{propNc}
The symmetric bilinear form $(\ | \ )_c$ is $W$-invariant. It is nondegenerate
on $V_c$ for all values of $m$ except a finite number of (absolutely real) algebraic
integers. If moreover $\k\subset \R$, it is a
scalar product, and in particular is nondegenerate, for all $m > N(c)$.
The discriminant of $(\ | \ )_c$ is a polynomial in $m$
of degree $\# c$, monic up to a sign, which admits $(m -N(c)$
as multiplicity 1 factor.    
\end{prop}
\begin{proof}
The fact that $(\ | \ )_c$ is $W$-invariant is an immediate
consequence of $\alpha(wsw^{-1},wuw^{-1}) = \alpha(s,u)$. Its
discriminant is a characteristic polynomial in $m$, hence
is monic up to a sign. Since the entries of $A_c$ are
integers, its roots are algebraic integers. The other assertions
are elementary consequences of the above lemma.
\end{proof}

\begin{remark}
\end{remark}
Calculations show (see appendix) that the roots should actually
be rationals (hence integers). We did not find a general argument
to justify that, though.

\medskip

\begin{prop} \label{propsP} For $s \in \RR$ we let $p_s \in \End(V)$ be defined by
$p_s.v_u = (v_s|v_u) v_s$.
\begin{enumerate}
\item $p_s^2 = (1-m) p_s$.
\item $p_s$ is selfadjoint with respect to $(\ | \ )$.
\item if $m \neq 1$ then $\frac{1}{1-m}p_s$ is the orthogonal projection
on $\k v_s$ with respect to $(\ | \ )$.
\item $sp_s = p_s s = p_s$.
\item for $w \in W$ we have $wp_sw^{-1} = p_{wsw^{-1}}$.
\end{enumerate}
\end{prop}
\begin{proof}
Since $(v_s|v_s) = 1-m$ then $p_s^2.v_u = p_s((v_s|v_u) v_s)
= (1-m) (v_s|v_u) v_s = (1-m) p_s.v_u$ for all $u \in \RR$,
which proves (1). Then
$(p_s.v_u|v_t) =  (v_s|v_u)(v_s|v_t) = (v_u | p_s.v_t)$
proves (2), and (3) is a consequence of (1) and (2). Then (4)
follows from the $W$-invariance of $(\ | \ )$ and, either a direct
calculation, or (3). Item (5) is clear.

\end{proof}

\section{Definition of the representation}

We define here a representation $\rho : \mathcal{T} \to \gl(V)$,
depending on $m \in \k$. In order to simplify notations,
for $s \in \mathcal{R}$ we
denote $t_s = t_H$ for $H = \Ker(s-1)$, and $t_s.x = \rho(t_s)(x)$
for $x \in V$.
We define $\rho$ as follows :
$$
t_s.v_s = m v_s, \ \ t_s.v_u = v_{sus} - \alpha(s,u) v_s \mbox{ for } s \neq u,
$$
where $\alpha$ is defined in the previous section. Using the
notations defined there, $t_s$ is defined to act like $s-p_s$.

\begin{prop} \label{omintegr} $\rho$ is a $W$-equivariant representation of $\mathcal{T}$
on $V$.
\end{prop}
\begin{proof}
By proposition \ref{propsP} it is clear that $w \rho(t_s) w^{-1} =
\rho(t_{wsw^{-1}})$, so we only have to check that the defining relations
of $\mathcal{T}$ are satisfied by the endomorphisms $\rho(t_s)$ for $s \in \RR$.
Let $Z = H_{s_1} \cap H_{s_2}$ be a codimension
2 subspace of $E$, and $\RR_0 = \RR_{s_1,s_2}$. We let
$t = \sum_{s \in \RR_0} t_s$. Let $x \in \RR_0$, and extend $\rho$to linear
combinations of $t_s,s \in \RR$. We
have to show that $[\rho(t),\rho(t_x)].v_u = 0$
for all $u \in \RR$.

We first prove that, if $u \in \RR_0$, then
$$
t.v_u = \left( m + C_0(u) -1\right) v_u \mbox{ where } C_0(u)
= \# \{ s \in \RR_0  \ | su=us \}.
$$
Indeed,
since $\alpha(s,u) = \# \{y \in \RR_0 \ | \ s = yuy \}$, we have
$$\sum_{s \in \RR_0 \setminus \{ u \} } v_{sus}
= \#\{ s \in \RR_0 \ | \ s\neq u \mbox{ and } sus = u \} v_u
+ \sum_{s \in \RR_0 \setminus \{ u \}} \alpha(s,u) v_s
$$
hence
$$
\begin{array}{lcl}
t.v_u &=& m v_u + \sum_{s \in \RR_0 \setminus \{ u \}} v_{sus} -
\alpha(s,u) v_s \\
 &= & m v_u + \sum_{s \in \RR_0 \setminus \{ u \}} v_{sus}
- \sum_{s \in \RR_0 \setminus \{ u \}}\alpha(s,u) v_s \\
&=& (m + \#\{ s \in \RR_0 \ | \ s\neq u \mbox{ and } sus = u \}) v_u \\
&=& (m + \#\{ s \in \RR_0 \ | \  sus = u \}-1) v_u \\
& = & (m + C_0(u) - 1) v_u
\end{array} 
$$
If $u=x$ it is clear that $[\rho(t),\rho(t_x)].v_u = 0$. Assuming $u \neq x$,
a direct computation then shows $$tt_x.v_u - t_x t .v_u = (C_0(xux) - C_0(u)) v_{xux}
-\alpha(x,u)(C_0(x) - C_0(u)) v_x.$$
On the other hand, for any $y \in \RR_0$, we have
$$C_0(yuy) = \# \{ s \in \RR_0 \ | \ syuys = yuy \} =
\# \{ s \in \RR_0 \ | \ ysyuysy = u \} = C_0(u).$$
In particular $C_0(xux) = C_0(u)$ and $[\rho(t),\rho(t_x)].v_u =
-\alpha(x,u)(C_0(x) - C_0(u)) v_x$. If $\alpha(x,u) = 0$
we are done, otherwise there exists $y \in \RR_0$ such
that $u = yxy$ and then $C_0(x) = C_0(yxy) = C_0(u)$.
It follows that $[\rho(t),\rho(t_x)].v_u = 0$ provided that $u \in \RR_0$.

We now assume $u \not\in \RR_0$. In particular $u \neq x$ and,
using $t.v_x = (m+C_0(x)-1)v_x$,
a direct calculation shows
$$
t t_x.v_u  = \sum_{s \in \RR_0} v_{sxuxs} - \sum_{s \in \RR_0 \setminus
\{ x \} } \alpha(s,xux) v_s - \alpha(x,xux) v_x
 - \alpha(x,u)(m+C_0(x) - 1) v_x
$$
Since $\alpha(x,xux) = \alpha(xxx,u) = \alpha(x,u)$ we have
$$
t t_x.v_u  = \sum_{s \in \RR_0} v_{sxuxs} - \sum_{s \in \RR_0 \setminus
\{ x \} } \alpha(s,xux) v_s 
 - \alpha(x,u)(m+C_0(x)) v_x
$$
Similarly, $t_xt.v_u$ equals
$$
 \sum_{s \in \RR_0} v_{xsusx} - \sum_{s \in \RR_0 \setminus
\{x \}} \alpha(s,u) v_{xsx}  - \left(m\alpha(x,u)+ \sum_{s \in \RR_0}
\alpha(x,sus)  - \sum_{s \in \RR_0 \setminus \{ x \}} \alpha(s,u)
\alpha(x,s) \right) v_x
$$
We claim that these two expressions are equal. First note that
$v_{xsusx} = v_{xsx(xux)(xsx)}$ and $s \mapsto xsx$ is a permutation
of $\RR_0$, hence
$$
\sum_{s \in \RR_0} v_{xsusx} = \sum_{s \in \RR_0} v_{s(xux)s}.
$$
Moreover, $\alpha(s,xux) = \alpha(xsx,u)$ implies
$$
\sum_{s \in \RR_0 \setminus \{ x \} } \alpha(s,xux) v_s = 
\sum_{s \in \RR_0 \setminus \{ x \} } \alpha(xsx,u) v_s
= \sum_{s \in \RR_0 \setminus \{ x \} } \alpha(s,u) v_{xsx}
$$
Finally, using $\alpha(x,sus) = \alpha(sxs,u)$, we get
$$
[t,t_x].v_u = \left( \sum_{s \in \RR_0} \alpha(sxs,u) -
\sum_{s \in \RR_0 \setminus \{ x \} } \alpha(s,u) \alpha(x,s) \ - \alpha(x,u) C_0(x) \right) v_x.
$$
We now consider the map $\varphi_x : \RR_0 \to \RR_0$
which maps $s$ to $sxs$. For $s' \in \RR_0$ with $s' \neq x$,
we have by definition $\# \varphi_x^{-1}(s') = \alpha(s',x)$,
and $\# \varphi_x^{-1}(x) = C_0(x)$. It follows that
$$
\sum_{s \in \RR_0} \alpha(sxs,u) = \sum_{s \in \RR_0\setminus \{x \} }
\alpha(s,x)\alpha(s,u) + C_0(x) \alpha(x,u)
$$
which implies that $[\rho(t),\rho(t_x)].v_u = 0$, namely that $\rho$ indeed defines
a representation
of $\mathcal{T}$.  
\end{proof}

It is readily checked that, for every $s \in \RR$,
any $V_c$ for $c \in \RR/W$ is setwise stabilized by $s,p_s$ and $t_s$.
We denote $\rho_c : \mathcal{T} \to \gl(V_c)$ the corresponding
representation. Whenever needed, we now identify $s,t_s$ for $s \in \RR$
with their action in the representation under investigation.
We obviously have
$$
V = \bigoplus_{c \in \RR/W} V_c, \ \ \ \ 
\rho = \bigoplus_{c \in \RR/W} \rho_c.
$$
Note that, if $m \neq -1$, for each $s \in \RR$
the $t_s$-stable subspaces of $V_c$
are $s$-stable. Indeed, $sp_s = p_s s$ and
$p_s^2 = (1-m) p_s$ by proposition \ref{propsP}, hence
from $t_s = s-p_s$ one easily gets
$$
s = \frac{-1}{m+1} t_s^2 + t_s + \frac{1}{m+1}.
$$

We prove the following.
\begin{prop} \label{propirredrhoc}
Let $c \in \RR/W$. Assume $m \neq -1$ and $(\ | \ )_c \neq 0$. Then $\rho_c$ is irreducible if and
only if $(\ | \ )_c$ is nondegenerate. In this case it is absolutely
irreducible.
\end{prop}

\begin{proof}
First assume that $(\ | \ )_c$ is degenerate, and let
$U = \Ker (\ | \ )_c \neq \{ 0 \}$. Since $(\ | \ )_c$ is
$W$-invariant, $U$ is setwise stabilized by $W$. Let
$s \in \RR$.
If $v \in U$ then $p_s .v = (v|v_s) v_s = 0$
and $U$ is also stabilized by $p_s$, hence by $t_s = s-p_s$.
It follows that $U$ is $\mathcal{T}$-stable. It is a proper
subspace since $(\ | \ )_c \neq 0$, hence $\rho_c$ is not
irreducible.

Conversely, we assume that $(\ | \ )_c$ is nondegenerate.
Without loss of generality we can assume that $\k$ is algebraically closed. Let $U \subset
V_c$ a $\mathcal{T}$-stable subspace with $U \neq \{ 0 \}$.
We show that $U = V$. Let $v \in U \setminus \{ 0 \}$.
Since $m \neq -1$ we know that $U$ is $W$-stable by the
above remark. Since $p_s = s-t_s$, it follows that $p_s .v = (v_s|v)
v_s \in U$ for all $s \in \RR$. If, for all $s \in c$, we had
$p_s.v = 0$, that is $(v_s|v)_c = 0$, then $v \in \Ker (\ | \ )_c = \{ 0 \}$,
a contradiction. Thus there exists $s_0 \in c$ such that
$(v_{s_0}|v)_c \neq 0$, hence $v_{s_0} \in U$. Since $U$ is
$W$-stable it follows that $v_s \in U$ for all $s \in c$
and $U = V$.

\end{proof}

\begin{remark}
\end{remark}
When $m = -1$, the same argument shows that the action of $\k W \ltimes
\mathsf{U} \mathcal{T}$ is irreducible iff $(\ | \ )_c$
is nondegenerate. A slight modification in the proof of
proposition \ref{propintegrirred} (or see the proof of \cite{REPTHEO}, propositions 7 or 8)
then shows that, when $\k = \C$,
the monodromy representation of $B$ is (absolutely) irreducible.
However, its restriction to $P$ need not be irreducible, as illustrates
the example of Coxeter type $A_2$. In this case,
$(\ | \ )$ is nondegenerate for $m = -1$ (see proposition \ref{propdiscri}) but
it is easily checked that $\rho$ admits two irreducible
components.

\bigskip

As a consequence of proposition \ref{propirredrhoc}, the sum of the $t_s$ for $s \in \mathcal{R}$, which
is central in $\mathcal{T}$, acts by a scalar on $V_c$ if $(\ | \ )_c$
is nondegenerate and $m \neq -1$. We show that this holds for all values of $m$,
regardless of the nondegeneracy of $(\ | \ )_c$.
Note that the cardinality of the centralizer $C_W(s)$ of $s$ in $W$ does
not depend on the choice of $s \in c$. We let $C(c) = \# C_W(s) \cap
\RR$ for some $s \in c$.

\begin{prop} \label{propTaction} The element $T = \sum_{s \in \RR} t_s \in \mathcal{T}$
acts on $V_c$ by the scalar $m -1+C(c)$.
\end{prop}
\begin{proof}
Let $u \in c$. We have
$$
T.v_u = m v_u + \sum_{s \in \RR\setminus \{ u \}} \left( v_{sus} - \alpha(s,u) v_s \right).
$$
We notice that $\alpha(s,u) = 0$ if $s \not\in c$. Moreover,
considering the map $\varphi_u : s \mapsto sus$ from $\RR$ to $c$,
we note that each $x \in c \setminus u$ has $\alpha(x,u)$ inverse images,
and $v_u$ has $\# C_W(s) \cap \RR - 1$ inverse images. It follows that
$$
T.v_u = m v_u + (C(c) -1) v_u +
\left( \sum_{x \in c \setminus \{ u \} } \alpha(x,u) v_x 
\right) - \left( \sum_{s \in c \setminus \{ u \} } \alpha(s,u) v_s 
\right) = (m + C(c) -1) v_u$$
and the conclusion follows.
\end{proof}

\begin{prop} \label{spectrets} Let $s \in \RR$, and assume that there exists $u \in \RR$
such that $su \neq us$ and $\alpha(s,u) \neq 0$. The action
of $t_s$ on $V$ is semisimple if and only if $m \neq 1$. In that case
it has eigenvalues $m,1,-1$, and we have $\Ker(s-1)=
\Ker(t_s-m)\oplus \Ker(t_s-1)$ and $\Ker(s+1) = \Ker(t_s+1)$.
\end{prop}
\begin{proof}
If $m = 1$, let $u \in \RR\setminus \{ s \}$ such that $sus \neq u$  and
consider the subspace $<v_s,v_u,v_{sus}>$. The matrix of $t_s$ on the
basis $(v_s,v_u,v_{sus})$ is
$$
\begin{pmatrix} 1 & -\alpha(s,u) & -\alpha(s,u) \\ 0 & 0 & 1 \\ 0 & 1 & 0
\end{pmatrix},
$$
which is not semisimple. It follows that the action of $t_s$ on $V$
is not semisimple either. Assume now $m \neq 1$. For $u \in \RR$,
we define $v'_s =v_s$ and $v'_u = v_u + \frac{\alpha(u,s)}{m-1} v_s$
for $u \neq s$. Then $t_s.v'_s = m v_s$ and
$$t_s.v'_u = v_{sus} - \alpha(u,s) v_s + \frac{\alpha(u,s)m}{m-1} v_s
= v_{sus} + \frac{\alpha(u,s)}{m-1} v_s = v'_{sus}.
$$
The $(v'_u)$ obviously form a basis of $V$. We have a partition $\RR = \{ s \} \sqcup \RR_0 \sqcup \RR_1$
with $\RR_0 = \{ u \in \RR \ | \ s\neq u, su = us \}$,
$\RR_1 = \{ u \in \RR \ | \ s\neq u, su \neq us \}$.
We have $V = \k v_s \oplus V_0 \oplus V_1$ where $V_i$ is
generated by the $v'_u$ for $u \in \RR_i$. Then $t_s$ acts by $m$
on $\k v_s$, $1$ on $\RR_0$ and $s$ acts by $1$ on $\k v_s \oplus V_0$.
The subspace of $V_1$ is a direct sum of planes $< v'_u,v'_{sus}>$
on which $t_s$ and $s$ both act by the matrix $\begin{pmatrix}
0 & 1 \\ 1 & 0 \end{pmatrix}$. The conclusion follows.
\end{proof}

\begin{remark}
\end{remark}
The irreducibility of $\rho$ on $V_c$ for generic $m$ implies the
irreducibility of its dual representation $\rho^*$. More precisely,
the action of $\mathcal{T}$ on $V^*$ with dual basis $(v_s^*)$ is given
by
$$
\left\lbrace \begin{array}{lcl}
t_s.v_s^* &=& m v_s^* - \sum_{u \in \RR \setminus \{s \} } \alpha(s,u) v_u^* \\
t_s. v_u^* &= & v_{sus}^* \mbox{ if } s \neq u
\end{array} \right).
$$
hence the span of $(v_s)^*$ for $s \in c$, naturally identified
with $V_c^*$, is setwise stabilized by $\mathcal{T}$, and is
acted upon by the dual action of $\mathcal{T}$ on $V_c$.
A consequence of proposition \ref{propirredtensor} below is that $\rho^*$
is not isomorphic to $\rho$, at least for generic $m$.

\bigskip

Recall that to any $W$-module
can be associated a representation of $\mathcal{T}$
where $t_s$ acts in the same way as $s \in \RR \subset W$
(see e.g. \cite{BMR} lemma 4.11). It seems that this remarkable
phenomenon has been first
noticed by I. Cherednik in special cases (see \cite{CHEREDNIK}).
For that reason, 
we call it the Cherednik representation associated
to this $W$-module.

Let $\WC$ be a maximal parabolic subgroup of $W$,
and $\RR_0$ its set of reflections.
We denote $\mathcal{T}_0$ the corresponding
holonomy Lie algebra.
Thanks to Steinberg theorem, we know
that $\RR_0 = \RR \cap \WC$. We let
$\rho_0$ denote the corresponding
representation on $V_0$. The
vector space $V_0$ is naturally
identified to a subspace of $V$, through
an inclusion of $\WC$-modules. Recall that the
holonomy Lie algebra $\mathcal{T}_0$ is identified
to a sub-Lie-algebra of $\mathcal{T}$ by lemma \ref{embedLie}.

\begin{prop} \label{restrKram}
The subspace $V_0$ is $\mathcal{T}_0$-stable,
and $\rho_{|V_0} = \rho_0$. Moreover,
for generic values of $m$, the restriction of $\rho_0$ to
$\mathcal{T}_0$ is
the direct sum of $\rho_0$ and of the Cherednik representation
associated to the permutation action of $\WC$ on $\RR \setminus \RR_0$.
\end{prop}
\begin{proof}
It is clear from the formulas that $V_0$ is $\mathcal{T}_0$-stable.
The fact that this action is the one we need amounts to
saying that, for $s,u \in \RR_0$ with $s \neq u$,
$$
\alpha(s,u) = \# \{ y \in \RR_0 \ | \ ysy=u \},
$$
meaning that, if $ysy=u$ for some $y \in \RR$,
then $y \in \RR_0 = \RR \cap \WC$.
Let $x \in E$ such that $\WC = \{ w \in W \ | \ w.x = x \}$.
Since $s,u \in \WC$ we have $x \in H_s \cap H_u$. We
proved that $ysy=u$ for some $y \in \RR$ implies $y \in \RR_{s,u}$,
hence $H_s \cap H_u$ is pointwise stabilized and in particular
$y.x = x$, namely $y \in \WC$.
If $(\ | \ )$ is nondegenerate, then $V_0$
admits an orthogonal subspace $U$.
For all $z \in U$ and $s \in \RR_0$, we have
$(z|v_s) = 0$ hence $p_s . z = 0$. It follows that
$t_s = s-p_s$ acts like $s$ on $U$, hence the
action of $\mathcal{T}$ on $U$ is a Cherednik
representation. Moreover, the images of
$v_s$ for $s \in \RR \setminus \RR_0$
in $V / V_0$ obviously form
a basis, on which $W$ acts by permutation.
The conclusion follows.
\end{proof}

\section{Decomposition of the tensor square}

We prove here that the alternating and symmetric squares
of $\rho_c$ are irreducible for generic values of $m$. As a preliminary,
we need to study the connections between the various
endomorphisms of $V \otimes V$ associated to a given
reflection $s \in \RR$. In order to simplify notations,
we simply denote $t_s$ the endomorphism $\rho(t_s)$, and $s$
the corresponding endomorphism of $V$. Recall that $t_s = s -p_s$.

\subsection{Endomorphisms of $V \otimes V$ associated to $s \in \RR$.}

Let $s \in \RR$. We let $T_s = t_s \otimes 1 + 1 \otimes t_s$,
$S = s \otimes s$, $\Delta_s = s \otimes 1 + 1 \otimes s$,
$P_s = p_s \otimes 1 + 1 \otimes p_s$, $Q_s = p_s \otimes p_s$,
$R_s = p_s \otimes s + s \otimes p_s$. We let $\mathcal{D}_s$
be the subalgebra with unit of $\End(V \otimes V)$ generated
by these six elements. Since $s,p_s,t_s,1$ commute with each other,
$\mathcal{D}_s$ is commutative. We have $T_s = \Delta_s - P_s$.

By direct computation, we checked that $\mathcal{D}_s$ has dimension
at most 6, with the following (commutative) multiplication table
$$
\begin{array}{c||c|c|c|c|c|}
 & \Delta_s & P_s & Q_s & R_s & S \\
\hline
\hline
\Delta_s & 2 + 2S & P_s +R_s & 2 Q_s & R_s + P_s & \Delta_s \\
 P_s & & (1-m)P_s+2Q_s & 2(1-m)Q_s & (1-m)R_s + 2Q_s & R_s \\
 Q_s & & & (1-m)^2 Q_s & 2(1-m) Q_s & Q_s \\
 R_s & & & & (1-m)P_s + 2 Q_s & P_s \\
 S & & & & & 1 \\
\hline
\end{array}
$$  
This enables us to express the powers of $T_s$ in terms of the other
endomorphisms. A linear algebra computation shows that
the elements $(T_s^k)$ for $0 \leq k \leq 5$ form
a basis of $\mathcal{D}_s$ as soon as $16m(m-3)(m+3)(m+1)^4 \neq 0$.
Explicitely, for every $m\in \k$, $m(m+3)(m-3)(m+1) P_s$
equals
$$
(25m^2-9)T_s - 30mT_s^2 + \frac{45-25m^2}{4} T_s^3 + \frac{15m}{2} T_s^4
- \frac{9}{4} T_s^5,
$$
{} $m(m+1)^2(m-3)(S+1)$
equals
$$(m+1)(5m+3)(m-1)T_s + \frac{1}{2} m(m^3-13m-19-m^2) T_s^2
- \frac{5m^3+3m^2-9m-15}{4} T_s^3 + m(m+2)T_s^4 - \frac{m+3}{4} T_s^5
$$
and finally
$$
m(m+1)^2 Q_s = \frac{1-m^2}{2} T_s + m T_s^2 + \frac{m^2-5}{8} T_s^3
- \frac{m}{4} T_s^4 + \frac{1}{8} T_s^5.
$$

\subsection{Endomorphisms associated to $(s,u) \in \RR^2$}

In order to simplify notations, for $X,Y \in \End(V)$
we let $X \bullet Y = X \otimes Y + Y \otimes X$. We let
$\mathcal{D}$ denote the subalgebra with unit of $\End(V \otimes V)$
generated by the $T_s, s \in \RR$. The previous
subsection showed that, for $m \not\in \{ -3,-1,0,3 \}$, the
endomorphisms $P_s = p_s \bullet 1$ and $Q_s = p_s \otimes p_s
\in \mathcal{D}_s \subset \mathcal{D}$.

In this paragraph, we will prove the following.
\begin{lemma} \label{lemsuD} Let $u,s \in \RR$ with $s \neq u$.
Then $p_s \bullet p_u \in \mathcal{D}$ provided that
$m\not\in \{ -3,-1,0,1,3 \}$
and $(m-1)^2 + 2 \alpha(s,u)^2 \neq 0$.
\end{lemma}
The following properties are easily checked :
\begin{enumerate}
\item $\forall s,u \in \RR\ \ p_up_sp_u = (v_s|v_u)^2 p_u$
\item $\forall s \in \RR \ \forall w \in W$, $s = wuw^{-1} \Rightarrow
p_sp_u = (v_s|v_u) w p_u.$
\item $\forall s,u \in \RR \ \ (p_sp_u)^2 = (v_s|v_u)^2 p_s p_u$.
\end{enumerate}
Let now $X = (p_s \bullet 1)(p_u \bullet 1) \in \mathcal{D}$.
We have $X = (p_s p_u)\bullet 1 + p_s \bullet p_u$. From
$p_x^2 = (1-m) p_x$ and (3) we get
$$
\begin{array}{lcl}
X^2 &= & (v_s|v_u)^2 p_sp_u \bullet 1 + (1-m)^2 p_s \bullet p_u +
2 p_sp_u \otimes p_s p_u + p_s p_up_s \bullet p_u + (1-m) p_sp_u \bullet 
p_u \\
&+ &(1-m) p_sp_u \bullet p_s + p_up_sp_u \bullet p_s + p_s p_u \bullet
p_up_s. \\
\end{array}
$$
From (1) we get
$$
\begin{array}{lcl}
X^2 &= & (v_s|v_u)^2 p_s p_u \bullet 1 + \left( (1-m)^2 + 2(v_s|v_u)^2 \right)
p_s \bullet p_u + 2 p_s p_u \otimes p_s p_u \\
&+& (1-m) p_s p_u \bullet p_u + (1-m) p_s p_u \bullet p_s +
p_sp_u \otimes p_up_s.\\
\end{array}
$$
Since $p_s \otimes p_s$ and $p_u \otimes p_u$ belong to $\mathcal{D}$
we have $(p_s \otimes p_s)(p_u \otimes p_u) = p_s p_u \otimes p_s p_u
\in \mathcal{D}$. Moreover $(p_s \bullet 1)(p_u \otimes p_u)
= p_s p_u \bullet p_u \in \mathcal{D}$ and
$(p_s \otimes p_s)(p_u \bullet 1) = p_s p_u \bullet p_s \in \mathcal{D}$.
It follows that
$$
Y = X^2 - (v_s|v_u)^2 X - 2 p_s p_u \otimes p_sp_u - (1-m)p_sp_u \bullet
p_u - (1-m) p_sp_u \bullet p_s \in \mathcal{D}
$$
and we just computed that
$$
Y = \left( (1-m)^2 + (v_s|v_u)^2 \right)p_s \bullet p_u + p_s p_u
\bullet p_u p_s.
$$
We assumed $s \neq u$. If
$\alpha(s,u) = 0$ and $s \neq u$, then $p_s p_u = p_s p_u = 0$,
hence $p_s \bullet p_u \in \mathcal{D}$ if $(1-m)^2 + \alpha(s,u)^2 \neq 0$.
Since $(1-m)^2 = (1-m)^2 + \alpha(s,u)^2 = (1-m)^2 + 2\alpha(s,u)^2$
this proves the lemma in that case.  
We now assume $\alpha(s,u) \neq 0$,
meaning that there exists $y \in \RR$ such that $yuy = s$. Since
$y = y^{-1}$ we get from (1) that
$p_s p_u = (v_s | v_u) y p_u$ and $p_u p_s = (v_s|v_u) y p_s$,
hence
$$
p_sp_u \bullet p_u p_s = (v_s|v_u)^2 (yp_u) \bullet (yp_s)
= (v_s|v_u)^2 (y \otimes y)(p_u \bullet p_s).
$$
It follows that
$$
\left(\begin{array}{c} Y \\ (y \otimes y)Y \end{array} \right)
= M
\left(\begin{array}{c} p_s \bullet p_u \\ (y\otimes y)p_s \bullet p_u \end{array} \right)
$$
with $M$ the following $2 \times 2$ matrix
$$
M = (1-m)^2 \mathrm{Id} + \alpha(s,u)^2 \begin{pmatrix} 1 & 1 \\ 1 & 1
\end{pmatrix}
$$
which has determinant $(1-m)^2( (m-1)^2 +2 \alpha(s,u)^2)$. This proves
that, under the assumption on $m$ given by the lemma, we have $p_s \bullet
p_u \in \mathcal{D}$, as stated.

\subsection{Irreducibility of $\Lambda^2 V_c$ and $S^2 V_c$}

We will use the following classical lemma.
\begin{lemma} \label{indform} Let $U$ be a finite-dimensional $\k$-vector space,
endowed with a nondegenerate symmetric bilinear form $<\ , \ >$.
The symmetric bilinear forms naturally induced by $<\ , \ >$
on $\Lambda^2 U$ and $S^2 U$ are nondegenerate.
\end{lemma}

This lemma can be proved for instance by first reducing to $\k$ algebraically closed,
then reducing to the standard form on $\k^n$, and finally explicitely
computing the matrix of the induced bilinear forms on a standard
basis ; or, by noting that they are restrictions of the induced bilinear form
on $U \otimes U$, which is clearly nondegenerate, to the (orthogonal)
eigenspaces of the selfadjoint operator $x \otimes y \mapsto y \otimes x$.

We will also need a graph-theoretic result. For a finite set $X$ and $r \geq 0$
we let $\mathcal{P}_r(X)$ denote the set of all its subsets of cardinality $r$.
If $\Gamma$ is a (nonoriented) graph on $X$, we associate to $\Gamma$ the
following graphs. We define $\Lambda^2 \Gamma$ to be a graph
with vertices $\mathcal{P}_2(X)$, with one edge between the elements
$\{ a,b \}$ and $\{a,c \}$ iff $b \neq c$ and there is an
edge in $\Gamma$ between $b$ and $c$. We define $S^2 \Gamma$
to be the graph with vertices $\mathcal{P}_1(X) \sqcup \mathcal{P}_2(X)$,
where $\{ a,b \} , \{a,c \}\in \mathcal{P}_2(X)$
are connected by an edge under the same condition as in $\Lambda^2 \Gamma$,
plus an edge between $\{ a \}$ and $\{a,b \} \in \mathcal{P}_2(X)$
iff there is an edge between $a$ and $b$ in $\Gamma$. We now prove
the following elementary result.

\begin{lemma} \label{SLconn} If $\Gamma$ is connected, then $\Lambda^2 \Gamma$
and $S^2 \Gamma$ are connected.
\end{lemma}
\begin{proof}
If $\# X \leq 2$ the statement is trivial,
hence we assume $\# X \geq 3$. Since $\Gamma$ is connected, for all $a \in X$
there exists $b \neq a$ which is connected to $a$ by an edge in $\Gamma$. It follows
that $\{ a \}$ is connected to some $\{a , b \}$ in $S^2 \Gamma$.
Since the restriction of $S^2 \Gamma$ to $\mathcal{P}_2(X)$ is
$\Lambda^2 \Gamma$, it is thus sufficient to show that $\Lambda^2 \Gamma$
is connected.

Let $\{a,b \}, \{a , b' \} \in \mathcal{P}_2(X)$ with $b \neq b'$.
We show that there exists a path in $\Lambda^2 \Gamma$
between $\{a,b \}$ and $\{a , b' \}$. If there exists a path
in $\Gamma$ between $b$ and $b'$ which does not pass through
$a$ then it provides a path between $\{a,b \}$ and $\{a , b' \}$.
Assuming otherwise, and because $\Gamma$ is connected,
we can choose a path $(b=b_0,b_1,\dots,b_n=b')$
in $\Gamma$ of shortest length between $b$ and $b'$. By assumption there exists
$r \in ]0,n [$ such that $b_r = a$. Since the path has shortest
length, $r$ is uniquely determined. Moreover,
$(b_r = a,b_{r+1},\dots,b_n = b')$ is a path in $\Gamma$ from $a$
to $b'$ which does not pass through $b$, since every path
between $b$ and $b'$ pass through $a$. This provides a
path in $\Lambda^2 \Gamma$ between $\{a,b \}$ and $\{b',b \}$.
Similarly, $(b = b_0,b_1,\dots,b_r = a)$ is a path in $\Gamma$ between $b$
and $a$ which does not pass through $b'$, hence a path in
$\Lambda^2 \Gamma$ between $\{b ,b'\}$ and $\{a , b' \}$.
It follows that, in all cases, there exists a path in $\Lambda^2 \Gamma$
between $\{a,b \}$ and $\{a , b' \}$.

Let now $\{a,b \}, \{c ,d \} \in \mathcal{P}_2(X)$ with 
$\{a,b \} \neq \{c ,d \}$. If $\# \{ a,b,c,d \} = 3$,
we proved that $\{a,b \}$ and $\{ c,d \}$ are connected
by a path in $\Lambda^2 \Gamma$. If $\# \{ a,b,c,d \} = 4$,
then there are paths between $\{a ,b \}, \{a,c \}$ and
$\{a,c \}, \{c, d \}$, hence between $\{a,b \}$ and $\{c,d \}$.
This proves that $\Lambda^2 \Gamma$ and $S^2 \Gamma$ are connected.

\end{proof}

We can now prove the main result of this section.
\begin{prop} \label{propirredtensor}
Except for a finite number of values $m$, the $\mathcal{T}$-modules
$\Lambda^2 V_c$ and $S^2 V_c$ are irreducible.
\end{prop}

\begin{proof}
Assume that $U \subset \Lambda^2 V_c$
is a $\mathcal{T}$-stable nonzero subspace, and $x \in U \setminus \{ 0 \}$.
We choose a total ordering on $c$. Then $x$ can be written as
$$
x = \sum_{s <u } \la_{s,u} v_s \wedge v_u
$$
where we denote $v_a \wedge v_b = v_a \otimes v_b - v_b \otimes v_a$.
For generic $m$ the form $(\ | \ )_c$ is nondegenerate, hence
the induced form on $\Lambda^2 V_c$
is nondegenerate by lemma \ref{indform}, and there exist $a,b \in c$ with $a<b$ such that
$(x|v_a \wedge v_b) \neq 0$. On the other hand, it is easily
checked that
$$p_a \bullet p_b (x) = \frac{1}{2} (x | v_a \wedge v_b)_c
v_a \wedge v_b
$$
hence $v_a \wedge v_b \in U$ by lemma \ref{lemsuD}, provided that $m$ do
not belong to the finite list of values excluded by this lemma. If $b' \in c \setminus \{ b,a \}$ with $\alpha(b',b) \neq 0$
we get 
$$p_a \bullet p_{b'} (v_a \wedge v_b)
= \left( (v_a |v_a)_c (v_b|v_{b'}) - (v_b|v_a)(v_a|v_{b'}) \right)
v_a \wedge v_{b'}
$$
Note that $(v_b|v_{b'})_c = \alpha(b,b') \neq 0$ by assumption, 
$(v_a | v_a)_c = 1-m$ and $(v_b|v_a)_c(v_a|v_{b'})_c = \alpha(a,b)\alpha(a,b')$.
It follows that $(v_a |v_a)_c (v_b|v_{b'})_c - (v_b|v_a)_c(v_a|v_{b'})_c \neq 0$
for generic values of $m$, hence $v_a \wedge v_{b'} \in U$.
We know that $\mathcal{G}_c$ is connected by lemma \ref{lemGconn}, hence
$\Lambda^2 \mathcal{G}_c$ is connected by lemma \ref{SLconn}.
It follows that $v_s \wedge v_{u} \in U$ for all $s \neq u$ in $c$.
This proves $U = \Lambda^2 V_c$, hence $\Lambda^2 V_c$
is irreducible for generic values of $m$. The proof
for $S^2 V_c$ is similar and left to the reader.
\end{proof}

\section{The monodromy representation of $B$}

\subsection{General result and main conjecture}

We let $\k = \C$, and denote $R$ (resp. $R_c$) the monodromy representation
of $B$ over $K$ associated to $\om_{\rho}$ (resp. $\om_{\rho_c}$)
and to an arbitrary basepoint $\underline{z}$.

By ``almost all $m$'' we mean for all but a finite number of
values of $m \in \C$. Recall that to every representation
of $W$ is naturally associated a representation of the
(cyclotomic) Hecke algebra of $W$, by monodromy of the corresponding
Cherednik representation (see \cite{BMR}). The unitary group $U_N^{\eps}(K)$
was defined in \S 2.

\begin{theo} \label{theorepmono}
The representation $R$ is isomorphic to the direct sum
of the $R_c$, for $c \in \RR/W$. We have $\dim R_c = \# c$. Moreover,
\begin{enumerate}
\item For almost
all $m$, and more specifically for $m > N(c)$, the
representation $R_c$ is absolutely irreducible. In that
case there exists a non-degenerate orthogonal form on $V_c$
for which $R_c(B) \subset U_{\# c}^{\eps}(K)$, where $V_c \otimes K$
is identified to $K^{\# c}$ through the basis $v_s, s \in c$. 
\item For almost
all $m$ the Zariski closure of $R_c(P)$ is $\GL(V_c \otimes K)$.
\item If $W_0$ is a maximal parabolic subgroup of $W$, then
for almost all $m$ 
the restriction of $R$ to $B_0$ is isomorphic to
a direct sum of the corresponding representation $R_0$
of $B_0$ and of the Hecke algebra representation
associated to the permutation representation of $W$
on $\RR \setminus \RR_0$.
\item Let $\sigma$ be a braided reflection associated to
some $s \in \RR$. We denote $q^{x} = \exp( x h)$, and
$K_c(s) = \# \{ u \in c \ | \ us \neq su \}/2$. Assume that there exists $u \in c$ such
that $su \neq us$ and $\alpha(s,u) \neq 0$, and that $m \neq 1$.
Then $R_c(\sigma)$ is semisimple with eigenvalues $q^m$ with multiplicity 1,
$-q^{-1}$ with multiplicity $K_c(s)$ and $q$ with multiplicity
$\# c -K_c(s)-1$.
\end{enumerate}
and we have
$$R_c(\bbeta)= q^{2 (m-1+C(c))/(\# c\ \# Z(W))}.$$
\end{theo}
\begin{proof}
Since $\rho$ is the direct sum of the $\rho_c$ for $c \in \RR/W$,
it is clear that $R$ is the direct sum of the $R_c$.
By proposition \ref{propirredrhoc} we know that $\rho_c$ is irreductible
as soon as $m \neq -1$ and $(\ | \ )_c$ is non-degenerate, and that in
this case it is absolutely irreducible. We know that this is in
particular the case for generic $m$ and when $m > N(c)$ by proposition \ref{propNc}.
In these cases $R_c$ is thus absolutely irreducible by
proposition \ref{propintegrirred}. We now prove the
assertion on unitarity :
we know that every $s \in \RR$ acts orthogonally
w.r.t. $(\ | \ )_c$ by proposition \ref{propNc},
hence is selfadjoint as $s^2 = 1$ ; it follows that $t_s = s-p_s$
is selfadjoint since $p_s$ is so by proposition \ref{propsP} (2),
so we can apply proposition \ref{propunitaire} to get the conclusion.

Now for almost all $m$ the $\mathcal{T}$-module $S^2 V_c$ and
$\Lambda^2 V_c$ are irreducible by proposition \ref{propirredtensor}.
Let $\g = \rho_c(\mathcal{T})$, for such an $m$.
If $\# c = 1$ we have $\dim V_c = 1$ and $\g' = \sl(V_c) = \{ 0 \}$.
We assume now $\# c > 1$. We have
$$
\End_{\g}(V_c\otimes V_c^*) \simeq (V_c \otimes V_c^*\otimes V_c^* \otimes
V_c )^{\g} \simeq  \End_{\g} (V_c \otimes V_c)
$$
which is 2-dimensional since $S^2 V_c \not\simeq \Lambda^2 V_c$
for dimension reasons, and $\# c \neq 1$. Moreover, since
$V_c$ is irreducible and faithful as a $\g$-module we know that $\g$ is
reductive. It follows that
the semisimple $\g$-module $\End(V_c)$ has two irreducible components,
hence $\sl(V_c)$ is irreducible as a $\g$-module.
Since $\g' = [\g,\g] \subset \sl(V_c)$ is non-zero we get
$\g' = \sl(V_c)$. By proposition \ref{propTaction} we know that
the sum $T$ of the $t_s$ acts on $V_c$ by the scalar $m - 1+C(c)$.
It follows that, for almost all $m$, $\g = \gl(V_c)$, hence
the Zariski-closure of $R_c(P)$ is $\GL(V_c \otimes K)$
by proposition \ref{propintegrirred}.

The assertion (3) on parabolic subgroups is a direct consequence of
theorem \ref{monorestr} and proposition \ref{restrKram}.
The assertion (4) is a consequence of proposition \ref{monobraidref} and
proposition \ref{spectrets}, and an easy computation of
$\dim \Ker(s+1)$ on $V_c$. Finally the computation of $R_c(\bbeta)$ follows from proposition \ref{monobeta}
and proposition \ref{propTaction}.
\end{proof}
If $W$ is a irreducible Coxeter group of type $ADE$ we showed in
\cite{KRAMMINF} that $R$ is isomorphic, after extension of scalars,
to a representation described in \cite{CGW} as the generalized
Krammer representation for types ADE, with
parameters $r = e^{- h}$, $l = e^{-m h}$. They are not
rigorously the
same than the ones described in \cite{DIGNE} and the original
Krammer one for type A, as it is described in \cite{KRAM}. The relation
is as follows. For an arbitrary representation $S$ of $B$
and $\la \in K$ we can define $\la S$ to be a representation
that sends every braided reflection $\sigma$ to $\la R(\sigma)$,
by using the morphism $\varphi : B \to \Z$ of remark \ref{remcentre}.
Then $e^{- h} R$ is isomorphic to the representation
described in \cite{DIGNE} and \cite{KRAM} for type A, with
parameters
$q = e^{-2 h}$ and $t = e^{(m+3) h}$.

In particular, this result obtained in \cite{KRAMMINF}
in conjunction with the faithfulness results of \cite{KRAM,CW,DIGNE}
can be translated as follows.

\begin{theo} If $W$ is a Coxeter group of type $ADE$ and $m \not\in \Q$,
then $R$ is faithful.
\end{theo}

This raises the question of whether the representations
obtained here are faithful in general. We believe that a positive answer
can be expected at least if $\# \RR/W=1$, whence conjecture 1 of the
introduction. More precisely, we conjecture the following stronger form.

\begin{conj} \label{mainconj}
If $W$ admits only one conjugacy class of reflections,
then $R$ is faithful for $m \not\in \Q$.
\end{conj}

This conjecture, if true, would imply group-theoretic properties
for arbitrary complex braid groups, that we investigate in section 8.
In section 7 we prove it for dihedral groups and give a direct
proof for Coxeter groups of type ADE, when $m$ is a formal
parameter. Note that, in the Coxeter
case, only the cases of types $H_3$ and $H_4$ remain open.

When $\# \RR/W > 1$, there may exist $c \in \RR/W$ such that $R_c$
is not faithful for any $m$. There is an obvious reason for this,
namely that, if $W$ is not irreducible, the representation
$R_c$ will factorise through the braid group associated
to the irreducible component of $W$ containing $c$. But even
if $W$ is irreducible this may happen, as we show below.

We do not have any guess concerning the following question.
A positive answer would prove the group-theoretic conjecture 2
of section 8. A negative answer may have its origin in the fact
that, when $\# \mathcal{R}/W $ is bigger than 1, the representation variety
of $B$ has usually larger dimension and thus faithfulness might be achieved
only for a more general representation depending on more parameters.

\begin{quest} If $W$ is irreducible, does there
exist $c \in \RR/W$ such that $R_c$ is faithful for $m \not\in \Q$~?
\end{quest}

\subsection{A nonfaithful component in type $B_n$}

Consider $W$ of type $B_n$
for large $n$. It has two classes of reflections. Let $c$ be its class of reflections of cardinality $n$,
and $x_1,\dots,x_n$ be the corresponding basis vectors of $V$.
Then $\mathcal{T}$ has generators $t_1,\dots,t_n$, $t_{ij}, t'_{i,j}$
for $1 \leq i,j \leq n$ and $i \neq j$, $t^{(')}_{ij} = t^{(')}_{ji}$, $t_{ii} = t'_{ii} = 0$.
The correspondance with the reflections of $W$ is as follows
$$
\begin{array}{lclcl}
t_{ij} & : & (z_1,\dots,z_i,\dots,z_j,\dots,z_n) & \mapsto & (z_1,\dots,z_j,\dots,z_i,\dots,z_n) \\
t'_{ij} & : & (z_1,\dots,z_i,\dots,z_j,\dots,z_n) & \mapsto & (z_1,\dots,-z_j,\dots,-z_i,\dots,z_n) \\
t_{i} & : & (z_1,\dots,z_i,\dots,z_n) & \mapsto & (z_1,\dots,-z_i,\dots,z_n) \\
\end{array}
$$
Then the action $\rho_c$ is given by
$$
\left\lbrace 
\begin{array}{lcll}
t_i.x_j & = & x_j - 2 x_i & \mbox{ if $i\neq j$} \\
t_i.x_i & = & mx_i\\
\end{array}
\right.
\left\lbrace 
\begin{array}{lclcll}
t'_{ij}.x_k & = & t_{ij}.x_k & = & x_k & \mbox{ if $k \not\in \{i,j \}$} \\
t'_{ij}.x_j & = & t_{ij}.x_j & = & x_i\\
\end{array} \right.
$$
Now consider the maximal parabolic subgroup of type $A_{n-1}$
that fixes the vector $(1,\dots,1)$. The restriction of $\rho_c$ to
$A_{n-1}$ is given by
$$\left\lbrace 
\begin{array}{lcll}
t_{ij}.x_k & = & x_k & \mbox{ if $k \not\in \{i,j \}$} \\
t_{ij}.x_j & = & x_i\\
\end{array} \right.
$$
which is the Cherednik representation associated to
the natural permutation action of $\mathfrak{S}_n$.
It is known that the corresponding representation
of the classical braid group on $n$ strands is the
(unreduced) Burau representation, and that this representation
is not faithful for large $n$ by work of Moody (see \cite{MOODY}).
By theorem \ref{monorestr} it follows
that $R_c$ is not faithful. On the other hand, both
representations $R_c$ and $R_{c'}$ are faithful on the center of $B$
by propositions \ref{monobeta} and \ref{propTaction}.
We showed in \cite{KRAMMINF} (see also section 8
below) that two normal subgroups of $B$ not included in its
center necessarily intersect each other. It follows that
$R$ is faithful if and only if $R_{c'}$ is faithful, where $\RR = c \cup c'$.

Note that the faithfulness of $R_{c'}$ would provide a faithful
representation of smaller dimension than the one deduced in \cite{DIGNE},
using so-called ``folding morphisms'', from the Krammer representation in type $A$,
and that the argument used here is already valid for $n \geq 5$ by \cite{BIGBUR}.

\subsection{Isomorphism class of the type $B_n$ components}

The two reflection classes $c,c'$ of $W$ in Coxeter type $B_n$
have cardinality $n$ and $n(n-1)$, respectively. We call the corresponding
irreducible components of $R$ the small and large component,
respectively. We let $\tau, \sigma_1,\dots,\sigma_{n-1}$
denote the standard Artin generators of $B$, with $<\sigma_1,\dots,
\sigma_{n-1}>$ the standard parabolic subgroup of type $A_{n-1}$
and $<\tau,\sigma_1>$ of type $B_2$.

\subsubsection{The small component}

We assume $n \geq 2$. The discriminant vanishes iff $m=-1$ or $m = 2n-1$ by proposition \ref{propdiscri}.
On this component, using the description of $\rho_c$ above and proposition \ref{spectrets}
we know that, if $m \neq 1$, then $\tau$ acts with eigenvalues $q, q^m$
and $\sigma_i$ with eigenvalues $q,-q^{-1}$, with both
actions being semisimple. Thus $R_c$ factors through the Hecke
algebra of type $B_n$ specialized at the corresponding (unequal) parameters.
The semisimplicity or genericity criterium for these
values of the parameter is satisfied as soon as
$|m-1| \geq 2n$ or $m-1 \not\in \Z$ (see e.g. \cite{ARIKI}).
Under these conditions (which imply $m \not\in \{1,-1\}$ hence $q^m \neq q$)
the representations of the Hecke algebra are
in 1-1 correspondance with couples
of partitions of total size $n$. The convention on this correspondance
depends on some ordering of the parameters. We choose the convention
such that the 1-dimensional representation
$([n],[0])$ is $\tau \mapsto q^m$,$\sigma_i \mapsto q$
and $([0],[n])$ is $\tau \mapsto q$,$\sigma_i \mapsto q$.
Now $\tau$ admits $q^m$ as eigenvalue with multiplicity 1, $q$ with multiplicity
$n-1$ and $\sigma_i$
admits $q$ as eigenvalue with multiplicity $n-1$. Since $R_c$ is irreducible
of dimension $n$,
according to the Hoefsmit models (see e.g. \cite{GP} \S 10.1) the only possibility is that
$R_c$ corresponds to the triple $([1],[n-1])$. We thus proved
the following.

\begin{prop} Assume that $W$ has type $B_n, n \geq 2$, that $|m - 1| \geq 2n$
or $m-1 \not\in \Z$,
that $m \neq 2n-1$, and that $c \in \mathcal{R}/W$ is
as above. Then the irreducible representation
$R_c$ factors trough the Hecke algebra of type $B_n$
with unequal parameters $q^m,q$ and $q,-q^{-1}$, and corresponds
to the couple of partitions $([1],[n-1])$.
\end{prop}

\subsubsection{The large component}

We assume $n \geq 2$. When $m \neq -1$, the
discriminant vanishes iff $m\in \{ 4n-5,2n-5 \}$, by proposition \ref{propdiscri}.
In \cite{HO}, H\"aring-Oldenburg introduced a family $BB_n(D)$
of finite-dimensional $D$-algebras, where $D$ is an integral commutative
$\C$-algebra with specified elements, used as parameters in the
definition of the algebra. In particular, there are distinguished units
$q,q_0,\la \in D$, where $q$ will ultimately be given the same meaning
as before.
When $q \neq q^{-1}$ these algebras are quotients of the group algebra $D B$
where $B$ is the Artin group of type $B_n$. The defining relations
for this quotient can be divided into the following list
\begin{enumerate}
\item Order relations. These are cubic relations $(\sigma_i - \la)(\sigma_i - q)(\sigma_i+q^{-1})=0$
and quadratic relations $\tau^2 = q_1 \tau + q_0$ with $q_1 \in D$.
\item Relations of type $A_{n-1}$. These are the additional relations
between $\sigma_1,\dots,\sigma_{n-1}$ which define the Birman-Wenzl-Murakami
algebra of type $A_{n-1}$.
\item Relations of type $B_2$. An additional relation involving
only $\tau$ and $\sigma_1$.
\end{enumerate}
In \cite{HO}, this algebra is shown to be semisimple and to have a nice structure
over the field of fractions $\tilde{D}$ of $D$
if
$q_0 = q^{-1}$, $1-q^{-1}\la$ is invertible in $D$, and there
exists a deformation morphism $D \to \C$ such that $q \mapsto 1, q_1 \mapsto 0$
(so called ``classical limit'').
More precisely,
its irreducible components are then in bijection with pairs $(a,b)$ of Young
diagrams whose sizes have for sum some integer at most $n$ and of the same
partity as $n$. Moreover, $BB_{n-1}(D)$
embeds in $BB_n(D)$ with multiplicity free restriction rule, which
is to add or remove one box in one of the Young diagrams.
Finally, it admits as quotient the Hecke algebra of type $B_n$
with relations $(\sigma-q)(\sigma+q^{-1}) = 0$ and $\tau^2  = q_1 \tau + q_0$.
The irreducible representations that factor through this quotient are indexed by the
couples of diagrams whose total size is equal to $n$.

Letting $D = \C[[h]]$, we choose for parameters $q = e^{ h}$, $q_0 = q^{-1}$,
$\la = q^m$ and $q_1 = q^{-1/2}(q-q^{-1})$, so that the order relation
on $\tau$ reads $(\tau - q^{1/2})(\tau +q^{-3/2})=0$. The deformation
map is given by $h \mapsto 0$. The condition that $1 - q^{-1}\la $
is invertible reads $m \neq 1$. Our convention
on the couples of partitions is that $([n],0)$ corresponds to
the 1-dimensional representation $\tau \mapsto q^{1/2}$
and $\sigma_i \mapsto q$.

We will identify, up to some renormalization, the representation
$R_{c'}$ with some component of this algebra. This should be compared
with the identification of $R$ in types $ADE$ with
an irreducible component of the generalized Birman-Wenzl-Murakami algebra (see \cite{KRAMMINF}).

\begin{prop} Assume $W$ has type $B_n, n \geq 2$, $m \not\in \{ -1,1 , 4n-5,2n-5 \}$, and $c' \in \mathcal{R}/W$
as above. Define $S(\sigma_i) = R_{c}(\sigma_i)$
and $S(\tau) = q^{-1/2} R_{c'}(\tau)$. Then the irreducible representation
$S$ factors through the H\"aring-Oldenburg algebra and corresponds to the
couple of partitions $([n-2],[0])$.
\end{prop}

The proof of this proposition is in the same spirit as in \cite{KRAMMINF},
section 4. Theorem \ref{theorepmono} states that $(R_{c'}(\tau)-q)(R_{c'}(\tau)+q^{-1})=0$
and $(R_{c'}(\sigma_i) - \la)(R_{c'}(\sigma_i)-q)(R_{c'}(\tau)+q^{-1})=0$,
whence $S(\tau)$ and $S(\sigma_i)$ satisfy the order relations (1).
In order to check the relations of type (2) we investigate the
restriction of $R_{c'}$ to type $A_{n-1}$. Let $v_{ij},v'_{ij}$
denote the natural basis elements of $V$, corresponding to the elements
$t_{ij},t'_{ij}$ of $\mathcal{T}$, and let $\mathcal{T}_0$
be the parabolic Lie subalgebra generated by the $t_{ij}$.
The subspace $U$ spanned by the $v_{ij}$ is stable by $\rho(\mathcal{T}_0)$
and is easily checked to be isomorphic to the infinitesimal Krammer representation of type $A_{n-1}$.
Now the action of $\mathcal{T}_0$ and $W_0$ on $V/U$ is readily seen
to be the Cherednik system associated to the permutation of the $v'_{ij}$.
Since the action of $\mathcal{T}_0$ on $V$ is semisimple, it follows from theorem \ref{monorestr} that
the restriction of $R_{c'}$ is the direct sum of the Krammer representation
of type $A_{n-1}$ and of a Hecke algebra representation. As a
consequence, it factorizes through the BMW algebra and relations
(2) are satisfied by $R_{c'}$ and $S$.

It remains to show that relations (3) are satisfied. This only depends
on the restriction of $R_{c'}$ to the parabolic subgroup of type $B_2$,
so in view of theorem \ref{monorestr} we investigate the corresponding restriction of $\rho_{c'}$
to the action of $t_1$ and $t_{12}$. Letting $w^{(')}_{ij} = v^{(')}_{ij} + \frac{1}{m-3} (v_{12} + v'_{12})$,
we have a direct sum decomposition of $V$ into the following subspaces,
on which $<t_{12},t_1>$ acts irreducibly :
\begin{itemize}
\item the plane $< v_{12}, v'_{12} >$
\item the planes $< w_{1j}- w'_{1j} , w_{2j} - w'_{2j} >$
\item the lines spanned by $w_{1j} + w'_{1j}+ w_{2j} + w'_{2j}$, $v_{ij}$, $v'_{ij}$
for $i,j \not\in \{ 1 ,2 \}$
\item the lines spanned by $
w_{1j} + w'_{1j}- w_{2j} + w'_{2j}$ for $j \not\in \{ 1 , 2 \}$
\end{itemize}
It is easily checked that the last three types actually factorize
through the Hecke algebra of type $B_2$, and correspond to
the pairs $([2],[0])$, $([1,1],[0])$ and $([1],[1])$.
The remaining representation factorizes through another Hecke algebra
of type $B_2$, namely the one with relation $(\sigma- q^m)(\sigma-q) = 0$
and the usual relation on $\tau$. Since this Hecke algebra
admits only one 2-dimensional irreducible representation, it is sufficient
to check that $BB_n$ admits such a representation. Recalling
from \cite{HO} the superfluous generator
$e_1 = 1 - \frac{1}{q-q^{-1}}(\sigma_1 - \sigma_1^{-1}) \in BB_n$
and the relations $\sigma_1 e_1 = \la e_1$, $\tau \sigma_1 \tau e_1 = e_1$
we get that $e_1,\tau e_1 \in BB_n$ span a subspace which is stable
by left multiplication. The nondegeneracy condition of \cite{HO}
implies that $e_1,\tau e_1$ are linearly independant.
We thus get a 2-dimensional representation of $BB_n$
whose matrix model on the basis $(e_1,\tau e_1)$ is given by
$$
\tau \mapsto \begin{pmatrix} 0 & q_0 \\ 1 & q_1 \end{pmatrix}
\ \ \ 
\sigma_1 \mapsto  \begin{pmatrix} \la & -q_1q_0^{-1} \\ 0 & q_0^{-1}
 \end{pmatrix}
$$
hence is irreducible and factors through this new Hecke algebra of type
$B_2$. This concludes the proof that the relations of $BB_n$
are satisfied.


We now prove that this irreducible representation of $BB_n$ corresponds to the
couple $([n-2],[0])$. For $n = 2$,
there is only one 2-dimensional irreducible representation of $BB_n$
which do not factor through the (usual) Hecke algebra, and it corresponds
to the couple $([0],[0]) = ([n-2],[0])$. 


We proceed by induction on $n$,
assuming $n \geq 3$. Then $S$ corresponds to some
couple of Young diagrams $(a,b)$ of total size at most $n$.
Since $S(\sigma)$ has 3 eigenvalues, $S$ does not factorize through the Hecke
algebra hence $|a|+|b| < n$. By theorem \ref{theorepmono}, the restriction of $S$ to
the parabolic subgroup of type $B_{n-1}$ contains the Hecke algebra
representation $([n-1],[0])$, and exactly one component which does not
factor through the Hecke algebra. The restriction rule for $BB_{n-1} \subset
BB_n$ thus implies $|a| + |b| = n-2$, $a \subset n-1$, $b \subset [0]$
hence $S$ corresponds to $([n-2],[0])$. This concludes the
proof of the proposition.

\subsection{Dihedral types}

We prove here the following result.

\begin{prop} \label{propfaithI2} If $W$ has type $G(e,e,2)$ with $e$ \emph{odd}, then
$R$ is faithful for generic $m$.
\end{prop}

We now assume that $W$ has type $G(e,e,2)$ with $e$ odd. It is known
by work of C. Squier (see \cite{SQUIER}) and has been reproved using another method
by G. Lehrer and N. Xi (see \cite{LEHXI}) that,
when $W$ is of dihedral type, then the suitably renormalized (reduced) Burau representation
of $B$ is faithful. Recall that the (reduced) Burau representation $R_b$
is the monodromy representation associated to the Cherednik system
on the reflection representation of $W$. The result is that
$R_b'$ defined by $R_b'(\sigma) = q R_b(\sigma)$ for $\sigma$
a braided reflection, is faithful. We use this result here,
by specializing the parameter $m$, in order to prove faithfulness of $R$.

We will actually prove the proposition when $m$ is a formal
parameter. This implies that $R$ is faithful for $m$ outside
a countable set of complex numbers. Indeed, for each $g \in B$
the matrix $R(g)$ has entries of the form $\sum_{k \geq 0} P_k(m) h^k$
where the $P_k$ are polynomials in $m$ (actually of degree at most $m$).
Each $g \in B \setminus \{ 1 \}$ thus defines a finite set of values
of $m$ for which $R(g) = \Id$. Since $B$ is countable it follows
that $R$ is faithful for $m \in \C$ outside the countable union of such
sets.

We now prove the proposition for $m$ a formal parameter.
First note that the naive approach to let $m = 1$, which
ensures that braided reflections have only two eigenvalues ($q$ and $-q^{-1}$), fails
here because in that case their images are not semisimple (by
proposition \ref{spectrets}).
In particular the specialisation of $R$ at $m = 1$ does not factor
through the Hecke algebra of $W$.

We use instead the specialization at $m = 0$, in which case
$\rho$ is not irreducible anymore by propositions \ref{propirredrhoc} and
\ref{propdiscri}, but admits as stable
submodule the kernel $U$ of $(\ | \ )$. By definition of $(\ | \ )$,
$t_s$ and $s$ act in the same way on $U$, for $s \in \mathcal{R}$.
It follows that the restriction of $\rho$ to $U$ is the Cherednik
representation associated to the action of $W$ on $U$, hence the
restriction of $R$ to $U \otimes K$ is the corresponding Hecke
algebra representation. It it thus sufficient to check that
$U$ contains a copy of the reflection representation of $W$.
Indeed, if this the case, then any $g \in \Ker R$ would lie
in $\Ker R_b$. Since $R'_b$ is faithful and $B$ is
generated by braided reflections then $R'_b(g) = q^r$ for
some $r \in \Z$, hence $g \in Z(B)$. Since $R$ is faithful
on $Z(B)$ for generic $m$ the conclusion would follow.

We now prove this representation-theoretic
property of the odd dihedral groups.
It is easily checked that $U = \{ \sum_{ s \in \mathcal{R}}
\la_s v_s \ | \ \sum \la_s = 0\}$. As $W$-modules, $V$ is the direct sum
of the trivial representation and $U$. It turns out that $U$ and $V$
are almost Gelfand models for $W$, proving that our assumption is valid.

\begin{lemma} As $\C W$-module, $U$ is the direct sum of
the 2-dimensional irreducible representations of $W$, all
occuring with multiplicity one.
\end{lemma}

\begin{proof}
Let $\chi_V$ and $\chi_U$ denote the characters of $V$ and $U$ as
$\C W$-modules, where $W$ has type $G(e,e,r)$, $e$ odd.
We have $\chi_U(g) = \chi_V(g)-1$ and $\chi_V(g) = \# \{ s \in \mathcal{R} \ | 
gs = sg \}$. It follows that, when $g \in W$, $\chi_U(g) = 0$
if $g$ is a reflection, $\chi_U(1) = e-1$ since $\# \mathcal{R} = e$,
and $\chi_U(g) = -1$ if $g$ is a nontrivial rotation.
Now $W$ admits $\frac{e-1}{2}$ irreducible 2-dimensional characters $\chi_k$,
for $1 \leq k \leq \frac{e-1}{2}$. They satisfy  $\chi_k(g) = 0$
if $g$ is a nontrivial rotation, $\chi_k(g) = \zeta^k + \zeta^{-k}$
if $g$ is a reflection and $\zeta$ some primitive
$e$-th root of 1. Since $e$ is odd, we have $\sum_k \zeta^k + \zeta^{-k} = -1$
hence $\chi_U = \sum_k \chi_k$ and the conclusion.
\end{proof}

This argument does not work in the case $\# \mathcal{R}/W = 2$,
as already shows the example of $G(4,4,2) = I_2(4)$. In that case,
the corresponding Hecke algebra representation has abelian image.

More generally, for types $G(2e,2e,2)$, this Hecke algebra
representation is not faithful, because it factorizes
through the (non-injective) morphism $B \to B'$, where $B'$ is the
braid group of type $G(e,e,2)$, which maps the Artin generators of $B$
to the Artin generators of $B'$. This can be seen as follows.
We denote by $W,W'$ the corresponding dihedral groups. The
Hecke algebra representations of $B$ which factor through a Hecke
algebra representation of $B'$ are, by Tits theorem,
precisely the deformations of the representations
of $W$ which factor through $W'$. Now the kernel of $W \to W'$ is
the center of $W$. Since this center acts trivially by conjugation on
the reflections, this proves that the representation of $W$ under
consideration indeed factors
through $W'$.

Note that $R$ itself may still be
faithful in these types. For instance it is faithful for type $I_2(4)$, as
it is easily checked to be the direct sum of two copies of the Burau representation
of type $B_2$. This argument however does not extend to other even
dihedral types.

\section{Group-theoretic properties}

We will show that conjecture \ref{mainconj}, if true, would imply a lot of properties for $B$,
beyond the usual consequences of linearity for a finitely
generated group (e.g. residual finiteness, Tits alternative, Hopf property, etc.).
It is usually difficult to prove the faithfulness of a representation
described as a monodromy of a local system. However, to get
the properties we have in mind, it would be sufficient
to construct a faithful representation with the same first-order
approximation as $R$, which may be easier. We make the statement precise in the
following theorem. Note that, for $g \in P$, the value of $R(g)$
modulo $h^2$ does not depend on the choice of the base point
$\underline{z} \in X$ ; indeed, $P$ is generated by loops $\gamma_H$ aroung the hyperplanes
$H \in \mathcal{A}$ (see \cite{BMR} prop. 2.2 (1))
and $R([\gamma_H]) \equiv \Id + 2h \rho(t_H)$ modulo $h^2$ (see the
proof of proposition \ref{propintegrirred}).

\begin{theo} \label{thprops}
Let $m \in \C$. Assume that $R$ is faithful or that
there exists a faithful representation $S$ of $B$ such that $S(g) \equiv R(g)$
modulo $h$ for all $g \in P$. Then
the following hold.
\begin{enumerate}
\item $P$ is residually torsion-free nilpotent.
\item $P$ is biorderable.
\end{enumerate}
Let $c \in \RR/W$ and $m \in \C$. Assume that $R_c$ is faithful or that there
exists a faithful representation $S$ of $B$ such that $S(g) \equiv R(g)$
modulo $h^2$ for all $g \in P$. Then the following hold.
\begin{enumerate} 
\item If $N_1, N_2$ are normal subgroups of $B$
such that $N_i \not\subset Z(B)$, then
$N_1 \cap N_2 \not\subset Z(B)$.
\item $B$ and its finite-index subgroup are
almost indecomposable in direct products,
meaning $G \simeq A \times B \Rightarrow A \subset Z(G) \mbox{ or }
B \subset Z(G)$.
\item The Fitting subgroup of $B$ (or $P$) equals its center.
\item The Frattini subgroup of $B$ (or $P$) is trivial, at least if $W \neq G_{31}$.
\end{enumerate}
\end{theo}

Before proving this result, we recall the group-theoretic notions involved
here. For some class $\mathcal{F}$ of groups, a group $\Gamma$
is called residually-$\mathcal{F}$ if, for all $g \in \Gamma \setminus \{ 1 \}$,
there exists $\pi : \Gamma \onto Q$ with $Q \in \mathcal{F}$ such that
$\pi(g) \neq 1$. The residual torsion-free nilpotence corresponds
to the class $\mathcal{F}$ of torsion-free nilpotent groups, and
is a significantly stronger property than torsion-free nilpotence.
For a finitely generated group it implies that the group $\Gamma$
is biorderable, namely that there exists a total ordering of the group
which is invariant by left and right multiplication (see \cite{PASSMAN} for
further consequences of this notion). It also implies that the  
group $\Gamma$ is residually a $p$-group for every prime $p$.

Moreover, this notion can be characterized in several ways.
Denote $C^r \Gamma$ the $r$-th term of the lower central
series of $\Gamma$. Being residually nilpotent means $\bigcap_r C^r \Gamma = \{ 1 \}$.
It is easily checked that being residually torsion-free
nilpotent means $\bigcap_r TC^r \Gamma = \{ 1 \}$ where $(TC^r \Gamma)$
is the torsion-free or rational lower central series, namely $TC^r \Gamma
= \{ g \in \Gamma \ | \exists n \in \Z \ | \  g^n \in C^r \Gamma \}$
is the preimage of the torsion subgroup in the nilpotent quotient $\Gamma
/C^r \Gamma$.

The Fitting subgroup of $\Gamma$ is the subgroup generated by all
normal nilpotent subgroups, and the Frattini subgroup is the
intersection of the maximal subgroups of $\Gamma$.

\begin{proof}
We have $R(P) \subset 1 + hM_N(A)$ for $N = \# \RR$ and $A = \C[[h]]$.
Since $1 + hM_N(A)$ is residually torsion-free nilpotent then (1) follows
if $R$ is faithful, or if $S(g) \equiv R(g)$ modulo $h^2$ for all
$g \in P$ and $S$ is faithful. Then residually torsion-free nilpotent
groups are biorderable, which implies (2).

Let $c \in \RR/W$. We showed that $R_c(P)$ is Zariski-dense in $\GL_N(K)$. The arguments
used to prove this apply verbatim to a representation $S$ as in
the statement of the theorem. It remains to show that, if $R : B \to
\GL_N(K)$ is faithful with $R(P)$ Zariski-dense in $\GL_N(K)$,
then the remaining properties (1) to (4) hold for $P$ and $B$. For
(1) to (3) this is proved in \cite{KRAMMINF}, theorem C and its first
corollary. For (4), corollary 2 of theorem C of \cite{KRAMMINF}
shows that the Frattini subgroup $\Phi(B)$ of $B$ is included in $Z(B)$.

Let $\varphi : B \to \Z$ be the morphism defined in remark \ref{remcentre},
which sends braided reflections to 1 and $\bbeta$ to
$2  \# \mathcal{R} /\# Z(W) \in \Z \setminus \{ 0 \}$.
For every prime number $p$, the group $H_p = \varphi^{-1}(p \Z)$
is clearly maximal, and $\cap_p H_p = \Ker \varphi$. It follows
that $\Phi(B) \subset \Ker \varphi$. If $W \neq G_{31}$ we know that $Z(B)$
is generated by $\bbeta$, hence $\Phi(B) \subset < \bbeta >\cap \Ker \varphi$,
which is trivial since $\varphi(\bbeta) \neq 0$. The same method
shows that $\Phi(P) = 1$, using $\bpi$ instead of $\bbeta$.

\end{proof}

The properties of $P$ and $B$ mentioned in this theorem
actually hold whenever $W$ is an irreducible Coxeter group,
not necessarily of type ADE, as proved in \cite{KRAMMINF}.
This suggests the following conjecture.

\begin{conj} \label{conjgroup} The groups $B$ and $P$ are
linear and they satisfy the properties enumerated in theorem \ref{thprops}
whenever $W$ is an irreducible
complex reflection group.
\end{conj}

Note that, if the properties of theorem \ref{thprops} hold
true for $W$ (say $W \neq G_{31}$) then the faithfulness of $R$ is equivalent
to the faithfulness of at least one of the $R_c$, $c \in \RR/W$.
Indeed, since the action of the center is faithful
through each $R_c$, the intersection $\Ker R$ of non-trivial kernels
$\Ker R_c$ cannot be trivial by property (1).

\bigskip

Conjecture \ref{conjgroup} is known for Artin groups by \cite{KRAMMINF},
that is when $W$ is a Coxeter group. It turns out that, when an irreducible $W$ has more
than one class of reflection, then $P$ and $B$ are closely connected to
(pure) Artin groups, at least when $\rk\, W \geq 3$. This enables to prove the
following, which provides additional support in favour of conjecture \ref{conjgroup}

\begin{theo} \label{thpropsB} The groups $B$ and $P$ are
linear and they satisfy the properties enumerated in theorem \ref{thprops}
whenever $W$ is an irreducible
complex reflection group with $\# \RR/W > 1$.
\end{theo}

This theorem is an immediate consequence of the two propositions below
and the corresponding results for Artin groups, using general
properties of Zariski-dense subgroups of $\GL_N(K)$ (see \cite{KRAMMINF},
\S 6.3).

\begin{prop} \label{propzarGL} If $W$ is irreducible and $\# \RR/W > 1$, then $B$ can
be embedded in some irreducible Artin group of finite Coxeter
type as a finite index subgroup. In particular, it can
be embedded in some $\GL_N(K)$ as a Zariski-dense subgroup.
\end{prop}
\begin{proof}
We use the classification of irreducible complex reflection groups.
The only exceptionnal types involved here are $G_{13}$ and $G_{28}$.
The type $G_{28}$ is the Coxeter type $F_4$.
It has been showed by Bannai \cite{BANNAI}
that the braid group $B$ of type
$G_{13}$ is isomorphic to the Artin group of type $I_2(6)$. Precisely,
$B$ has a
presentation $<x,y,z \ | \ yzxy=zxyz, zxyzx=xyzxy >$,
and the formulas $a = zx, b = zxy(zx)^{-1}$, and their inverse
$x = (baba)^{-1}, y = a^{-1}ba, z = (aba)^{-1}b(aba)$ provide
an isomorphism with the dihedral Artin group of presentation $<a,b \ | \ ababab=bababa>$
(these formulas were communicated to me by M. Picantin).

We now consider the infinite series. The types $G(e,e,2)$ with $e$ even
are (dihedral) Coxeter groups.
We are left with the $G(2e,e,r)$. In that case, $B$ can be embedded
in the Artin group of type $B_r$ as a finite index subgroup
of index $e$ (see \cite{BMR} \S 3 B1).
\end{proof}

Taking for $W$ a \emph{pseudo-}reflection group does not
enrich the collection of possible groups $B$, however it provides
new pure braid groups $P$. In that case we denote $\RR$ the set
of pseudo-reflections of $W$. Note, that if $W$ is not
a reflection group, we necessarily have $\# \RR/W > 1$.

We will prove the following.

\begin{prop} \label{propresnil} If $W$ is an irreducible pseudo-reflection group with $\# \RR/W > 1$,
then $P$ is residually torsion-free nilpotent.
\end{prop}

The proof uses the classification and arguments of several kinds.
The first argument is that, when $W$ has rank 2, then the
complement $X$ of the hyperplane arrangement is fiber-type
in the sense of Falk and Randell (see \cite{FALKRANDELL}),
hence $P = \pi_1(X)$ is residually torsion-free nilpotent
by \cite{FALKRANDELL,FALKRANDELL2}.

We thus can assume that $W$ has rank at least 3. Consider then
the infinite series $G(de,e,r)$ of complex pseudo-reflection
groups for $r \geq 3$. The assumption $\# \RR/W > 1$ is equivalent
to $d > 1$. But then the hyperplane arrangement is the same as the
one of type $G(de,1,r)$, which is fiber-type.

The remaining groups are the Artin group of type $F_4$,
for which the conclusion is known by \cite{KRAMMINF},
and the so-called Shephard group of rank 3, of type $G_{25}, G_{26}$
and $G_{32}$. For these we have to use another argument,
since none of them correspond to fiber-type arrangement
--- recall that the fiber-type condition is actually a combinatorial
one (see \cite{ORLIKTERAO}) and can be easily checked on a specific
lattice of hyperplanes. The Coxeter diagrams of these groups are the following ones.

\def\nnode#1{{\kern -0.6pt\mathop\bigcirc\limits_{#1}\kern -1pt}}
\def\ncnode#1#2{{\kern -0.4pt\mathop\bigcirc\limits_{#2}\kern-8.6pt{\scriptstyle#1}\kern 2.3pt}}
\def\sbar#1pt{{\vrule width#1pt height3pt depth-2pt}}
\def\dbar#1pt{{\rlap{\vrule width#1pt height2pt depth-1pt} 
                 \vrule width#1pt height4pt depth-3pt}}

$$
G_{25}\ \ \  \ncnode3s\sbar16pt\ncnode3t\sbar16pt\ncnode3u\ \ \ \ \ \ \ \ \ \ \ \ 
G_{26}\ \ \  \ncnode2s\dbar16pt\ncnode3t\sbar16pt\ncnode3u \ \ \ \ \ \ \ \ \ \ \ \ 
G_{32}\ \ \  \ncnode3s\sbar10pt\ncnode3t\sbar10pt\ncnode3u\sbar10pt\ncnode3v
$$

It is known (see \cite{BMR}) that removing the conditions on the order
of the generators gives a presentation of the corresponding braid group.
In particular, these have for braid groups the Artin groups of
Coxeter type $A_3, B_3$ and $A_4$, respectively.

We recall a matrix expression of the Krammer representation
for $B$ of Coxeter type $A_{n-1}$, namely for the classical braid group on $n$
strands. Letting $\sigma_1,\dots,\sigma_{n-1}$ denote its Artin generators
with relations $\sigma_i \sigma_j = \sigma_j \sigma_i$ if $|j - i| \geq 2$,
$\sigma_i \sigma_{i+1} \sigma_i = \sigma_{i+1} \sigma_i \sigma_{i+1}$, their action
on a specific basis $x_{ij}$ ($1 \leq i < j \leq n$) is given by (see \cite{KRAM})
$$
\left\lbrace \begin{array}{ll}
\sigma_k x_{k,k+1} = tq^2 x_{k,k+1} \\
\sigma_k x_{i,k} = (1-q)x_{i,k} + q x_{i,k+1} & i<k\\
\sigma_k x_{i,k+1} = x_{i,k} + t q^{k-i+1} (q-1) x_{k,k+1} & i<k \\
\sigma_k x_{k,j} = tq(q-1) x_{k,k+1} + q x_{k+1,j}& k+1<j \\
\sigma_k x_{k+1,j} = x_{k,j} + (1-q) x_{k+1,j} & k+1<j\\
\sigma_k x_{i,j} = x_{i,j} & i<j<k \mbox{ or } k+1<i<j \\
\sigma_k x_{i,j} = x_{i,j} + tq^{k-i} (q-1)^2  x_{k,k+1} & i<k<k+1<j\\
\end{array} \right.
$$
where $t$ and $q$ denote algebraically independent parameters. We
embed the field $\Q(q,t)$ of rational fractions in $q,t$ into $K = \C((h))$ by $q \mapsto -j e^h$ and
$t \mapsto e^{\sqrt{2}h}$, where $j$ denotes a primitive 3-root of 1.
We then check by an easy calculation that $\sigma_k^3 \equiv 1$
modulo $h$. Since the quotients of the braid group on $n$
strands by the relations $\sigma_k^3 = 1$
are, for $n = 3,4,5$, the Shephard group of types $G_4, G_{25}$ and $G_{32}$,
respectively, it follows that the pure braid groups of these types
embed in a residually torsion-free nilpotent group of the form $1 + h Mat_N(A)$ (where
$N = n(n-1)/2$) which proves their residual torsion-free nilpotence.

We are then left with type $G_{26}$. Types $G_{25}$ and $G_{26}$
are symmetry groups of regular complex polytopes which are
known to be closely connected (for instance they both appear in the
study of the Hessian configuration, see e.g. \cite{COXETER} \S 12.4 and
\cite{ORLIKTERAO} example 6.30). The hyperplane arrangement
of type $G_{26}$ contains the 12 hyperplanes of type $G_{25}$
plus 9 additional ones. The natural inclusion induces morphisms
between the corresponding pure braid groups, which cannot
be injective, since a loop around one of the extra
hyperplanes is non trivial in type $G_{26}$ (e.g. because $H^1(X)$
is freely generated by the forms $\om_H$ for $H \in \mathcal{A}$)
but becomes of course trivial in type $G_{25}$. 
However we will prove the
following, which concludes the proof of proposition \ref{propresnil}.

\begin{prop} \label{prop2526}The pure braid group of type $G_{26}$ embeds into
the pure braid group of type $G_{25}$.
\end{prop}

More precisely, letting $B_i,P_i,W_i$ denote the braid group, pure
braid group and pseudo-reflection group of type $G_i$, respectively, we
construct morphisms $B_{26} \into B_{25}$ and $W_{26} \onto W_{25}$
such that the following diagram commutes, where the vertical arrows
are the natural projections.
$$
\xymatrix{
B_{25} \ar@{->>}[d] & B_{26}\ar@{_{(}->}[l] \ar@{->>}[d]\\
W_{25} & \ar@{->>}[l]W_{26}
}
$$
Both horizontal morphisms are given by the formula
$(s,t,u) \mapsto ((tu)^3,s,t)$, where $s,t,u$ denote the generators
of the corresponding groups according to the above diagrams. The morphism between
the pseudo-reflection groups is surjective because it is a retraction
of an embedding $W_{25} \into W_{26}$
mapping $(s,t,u)$ to $(t,u,t^{sut^{-1}u})$. 
The kernel of this projection is the subgroup of order
2 in the center of $W_{26}$ (which has order 6).

We now consider the morphism between braid groups and prove
that it is injective. First recall that the braid group
of type $G_{26}$ can be identified with the Artin group of
type $B_3$. On the other hand, Artin groups of type $B_n$
are isomorphic to the semidirect product of the Artin group
of type $A_{n-1}$, that we denote $\mathcal{B}_n$
to avoid confusions,
with a free group $F_n$ on $n$ generators $g_1,\dots,g_n$,
where the action (so-called `Artin action') is given (on the left) by
$$\sigma_i :
\left\lbrace 
\begin{array}{lcll}
g_i &\mapsto &g_{i+1}\\
g_{i+1} &\mapsto & g_{i+1}^{-1} g_i g_{i+1} \\
g_j &\mapsto & g_j & \mbox{ si } j \not\in \{i,i+1 \}
\end{array} \right.  
$$
If $\tau, \sigma_1,\dots,\sigma_{n-1}$ are the standard generators of the Artin
group of type $B_n$, with $\tau \sigma_1 \tau \sigma_1 = \sigma_1 \tau \sigma_1 \tau$,
$\tau \sigma_i = \s_i \tau$ for $i > 1$, and usual braid relations
between the $\sigma_i$, then this isomorphism is given by
$\tau \mapsto g_1$,
$\sigma_i \mapsto \sigma_i$ (see \cite{CRISPPARIS} prop. 2.1 (2) for more details).
Finally, there exists an embedding of this semidirect product
into the Artin group $\mathcal{B}_{n+1}$ of type $A_{n}$ which satisfies $g_1 \mapsto
 (\sigma_2 \dots \sigma_{n})^n$, and $\sigma_i \mapsto \ \sigma_i$
($i \leq n-2$). By composing both, we get an embedding which makes
the square commute. This proves proposition \ref{prop2526}.

This embedding of type $B_n$ into type $A_n$,
different from the more standard
one $\tau \mapsto \sigma_1^2, \sigma_i \mapsto \sigma_{i+1}$,
has been considered in \cite{LONG}. The algebraic proof given
there being somewhat sketchy, we provide the details here. This
embedding comes from the following construction.

Consider the (faithful)
Artin action as a morphism $\mathcal{B}_{n+1} \to
\Aut(F_{n+1})$, and the free subgroup $F_{n} = < g_1,\dots,g_{n}>$ of $F_{n+1}$.
The action of $\mathcal{B}_{n+1}$
preserves the product $g_1 g_2\dots g_{n+1}$, and there is a natural
retraction $F_{n+1} \onto F_{n}$ which sends $g_{n+1}$ to $(g_1 \dots g_{n})^{-1}$.
This induces a map to $\Aut(F_{n})$, whose kernel is the center of
$\mathcal{B}_{n+1}$
by a theorem of Magnus (see \cite{MAGNUS}). Its image contains the group of inner automorphisms
of $F_{n}$, which is naturally isomorphic to $F_{n}$.

Indeed, it is straightforward to check
that $b_1 = (\sigma_2 \dots \sigma_{n})^{n}$
is mapped to $\Ad(g_1) = x \mapsto g_1 x g_1^{-1}$. Defining
$b_{i+1} = \sigma_i b_i \sigma_i^{-1}$, we get that
$b_i$ is mapped to $\Ad(g_i)$. In particular the subgroup
$\mathcal{F}_{n} = <b_1,\dots,b_{n}>$ of $\mathcal{B}_{n+1}$ is free and there is a natural
isomorphism $\varphi : b_i \mapsto g_i$ to $F_{n}$. Now,
let $\mathcal{B}_n \subset \mathcal{B}_{n+1}$ be generated by
$\sigma_i, i \leq n-1$.
For $\sigma \in \mathcal{B}_n$ and $b \in \mathcal{F}_{n}$ we know
that $\sigma b \sigma^{-1}$ is mapped to $\Ad(\sigma.\varphi(b))$
in $\Aut(F_{n})$, hence $\sigma b \sigma^{-1}$ and
$\varphi^{-1}(\sigma.\varphi(b)) \in \mathcal{F}_{n}$ may differ only by
an element of the center $Z(\mathcal{B}_{n+1})$ of $\mathcal{B}_{n+1}$.
On the other hand, $\varphi : \mathcal{F}_n \to F_n$
commutes with the maps $F_n \to \Z$
and $\eta : \mathcal{F}_{n} \to \Z$ which map every generator to 1.
Likewise, the
Artin action commutes with $F_n \to \Z$ hence
$\eta(\varphi^{-1}(\sigma.\varphi(b))) = \eta(b)$.
We denote $\ell : \mathcal{B}_{n+1} \to \Z$ the abelianization map. We have
$\ell(b_i) = n(n+1)$ for all $i$, hence $\ell(b) = n(n+1)\eta(b)$
for all $b \in \mathcal{F}_{n-1}$. Since $\ell(b) =
\ell(\sigma b \sigma^{-1})$
it follows that $\sigma b \sigma^{-1}$ and
$\varphi^{-1}(\sigma.\varphi(b)) \in \mathcal{F}_{n}$ differ by an element
in $Z(\mathcal{B}_{n+1}) \cap (\mathcal{B}_{n+1},\mathcal{B}_{n+1})$,
where $(\mathcal{B}_{n+1},\mathcal{B}_{n+1})$ denotes the commutators subgroup.
But $Z(\mathcal{B}_{n+1})$
is generated by $(\sigma_1 \dots \sigma_{n})^{n+1}
\not\in (\mathcal{B}_{n+1},\mathcal{B}_{n+1})$
hence $\sigma b \sigma^{-1}=\varphi^{-1}(\sigma.\varphi(b)) \in \mathcal{F}_{n}$.
In particular $\mathcal{F}_{n}$ is stable under the action by
conjugation of $\mathcal{B}_{n}$, which
coincide with the Artin action. This is the embedding $\mathcal{B}_n \ltimes \mathcal{F}_n \into
\mathcal{B}_{n+1}$ that is needed
to make the square commute.

\begin{remark} \end{remark}
We make an historical remark summarizing previously known results on
these topics. We first consider the residual torsion-free nilpotence of
$P$. After the seminal study of fiber-type arrangements
by Falk and Randell (see \cite{FALKRANDELL2}) and partial subsequent work
on Coxeter type $D$ (see \cite{MARK}) the case of
Coxeter arrangements has been settled in \cite{resnil,KRAMMINF}
by composing the idea that monodromy representations of this type
have the intersection of the lower central series in their
kernel and that the Krammer representation should behave like
(and actually is) a monodromy representation.
This same idea and a similarly short proof appeared
simultaneously (early 2005) in the work of
V. Leksin (see \cite{LEKSIN}, theorem 3;
note that $T C^r \Gamma$ coincides with the radical of $C^r \Gamma$,
where the radical of $H < \Gamma$ is the subgroup of $\Gamma$
generated by elements some power of which lies in $H$).
In \cite{LEKSIN} it is already mentionned that no additional work
is needed for the types $G(de,e,n)$ when $d > 1$ (the argument used
there, borrowed from \cite{BMR}, is not valid in the case $d=1$).

Concerning the useful property that normal subgroups
of $B$ which are not included in the center intersect each other, it
has been proved by D. Long in \cite{LONGNORMAL} for Artin groups of type $A$,
using their interpretation as mapping class groups of a punctured disk
and the Nielsen-Thurston classification of diffeomorphisms. Its
extension to irreducible Artin groups of finite Coxeter types was done in \cite{KRAMMINF}.

\section{Appendix : Computation of discriminants}

In this appendix, we compute the discriminants of $(\ | \ )$
for some irreducible groups of special interest (one discriminant for
each class of reflections).
In tables
\ref{table1} and \ref{table2}
we computed the discriminants for all exceptional groups,
as well as for irreducible groups of small size
in the infinite series $G(e,e,r)$ and $G(2e,e,r)$.
We were not able to prove a general formula for these groups, although
it is likely that there is an elementary one. For instance, the following
seems to hold.
\begin{conj}
If $e$ is odd, the discriminant for (the only reflection class of)
$G(e,e,r)$ is
$$
\pm (m-(2r-3)e)(m-(r-3)e)^{r-1}m^{(e-1)\frac{r(r-1)}{2}}
(m+e)^{\frac{r(r-3)}{2}}.
$$
\end{conj}

The next proposition gives the discriminant $\Delta_c$
of $(\ | \ )_c$ for the infinite series of Coxeter groups.

\begin{prop} \label{propdiscri} {\ \ }\\
\begin{enumerate}
\item If $W$ is a dihedral group of type $G(e,e,2)$, $e$ odd, then $\Delta = (-1)^e(m-e)m^{e-1}$.
\item If $W$ is a dihedral group of type $G(e,e,2)$, $e$ even, then $\Delta_c = (-1)^{\frac{e}{2}}(m-e+1)(m+1)^{\frac{e}{2}-1}$
for each $c \in \RR/W$.
\item If $W$ is a Coxeter group of type $A_{n-1}$, then
$$\Delta = (m+1)^{\frac{n(n-3)}{2}} (m-n+3)^{n-1}(m-2n+3).$$
\item If $W$ is a Coxeter group of type $G(2,2,n) = D_{n}$, then
$$\Delta = (m-4n+7) (m-1)^{\frac{n(n-1)}{2}} (m+3)^{\frac{n(n-3)}{2}}(m-2n+7)^{n-1}.$$
\item If $W$ is a Coxeter group of type $G(2,1,n) = B_{n}$, then
$$\Delta_{c_1} = (-1)^n(m-2n+1)(m+1)^{n-1} , \Delta_{c_2} = (m-4n+5)(m-2n+5)^{n-1}(m+1)^{n(n-2)}.
$$
\end{enumerate} 
\end{prop}
\begin{proof}
We use the presentation of the dihedral group $G(e,e,2)$ as $<s, \om \ | \ s^2 = \om^e = A, s \om = \om^{-1} s >$.
Then the reflections are the $s \om^i$, for $0 \leq i < e$. Now $s \om^i$ and $s \om^j$ are conjugated
by a reflection $s \om^a$ if and only if $i + j \equiv 2a$ modulo $e$. If $e$ is odd there
always exists such a $a$, which is then unique
modulo $e$. If $e$ is even there are exactly two of them if $s \om^i$ and $s \om^j$
are in the same conjugacy class.
It follows that, if $e$ is odd, the matrix of $(\ | \ )$ has
all off-diagonal entries equal
to 1. Likewise, if $e$ is even there are two classes of cardinality
$\frac{e}{2}$ and the matrix of each $(\ | \ )_c$ has all off-diagonal
entries equal to 2.
The determinants of such matrices being easy to compute, this proves (1) and (2).
(3) and (4) are proved in \cite{KRAMMINF} propositions 7.1 and 8.1, respectively. 
For (5), let us first consider the class $c_1$ formed by the
reflections
$$
r_i :
(z_1,\dots,z_i,\dots,z_n)  \mapsto  (z_1,\dots,-z_i,\dots,z_n) .
$$
Then $(v_{r_i} | v_{r_j}) = 2$
for $i \neq j$ and we conclude as in (1) or (2). Let now $c = c_2$.
We let $v_{ij}$ and $v'_{ij}$ be the basis vectors of $V_c$
corresponding to the reflections
$$
\begin{array}{lclcl}
v_{ij} & : & (z_1,\dots,z_i,\dots,z_j,\dots,z_n) & \mapsto & (z_1,\dots,z_j,\dots,z_i,\dots,z_n) \\
v'_{ij} & : & (z_1,\dots,z_i,\dots,z_j,\dots,z_n) & \mapsto & (z_1,\dots,-z_j,\dots,-z_i,\dots,z_n) \\
\end{array}
$$
It is easily checked that, if $i,j,k,l$ are distincts indices,
$(v_{ij}|v_{kl}) = (v_{ij}|v'_{kl}) = 0$ , $(v_{ij}|v_{jk})= (v_{ij}|v'_{jk})=1$
and $(v_{ij}|v'_{ij}) = 2$.
We let $V^0$ denote the subspace of $V_c$ spanned by the $v_{ij},v'_{ij}$ for $1 \leq i,j \leq n-1$,
$w_k = v'_{kn}-v_{kn}$, $y_k = v'_{kn}+v_{kn}$ for $k \leq n-1$. It is clear that the
$w_k,y_k$ span a supplement of $V^0$ in $V_c$. We check that $w_k \in (V^0)^{\perp}$,
$(w_k|w_l) = -2(m+1)\delta_{k,l}$, and $(w_k|y_l) = 0$ for all $k,l$. We introduce the
following elements of $V_0$ :
$$
u_k = \sum_{\{i,j \} \subset [1,n-1]} v'_{ij} + v_{ij},\ \ \  u = \sum_{1 \leq i,j \leq n-1} v'_{ij} + v_{ij}.
$$
Then one can check through an easy though tedious calculation that
$$
z_k = y_k - \frac{2}{2n-7-m} u_k + \frac{8}{(2n-7-m)(4n-9-m)} u \in (V^0)^{\perp}
$$
thus the $w_k,z_k$ span the orthogonal of $V^0$. Moreover $(w_k|z_l) = 0$
for all $k,l$ and $(z_k|z_l)$ only depend on whether $k=l$. It it thus easy to
compute the discriminant of $(\ | \ )$ on $(V^0)^{\perp}$. Up to some
non-zero scalar square, we find
$$
(m+1)^{2n-3} \frac{(m-2n+5)^{n-1}(m-4n+5)}{(m-2n+7)^{n-2}(m-4n+9)}
$$
and the conclusion follows easily by induction on $n$, since
we know that this discriminant is a polynomial in $m$ with
leading coefficient in $\{1,-1 \}$.
\end{proof}

Note that, for Coxeter groups of type ADE, we always have $\alpha(s,u) \leq 1$.
The situation changes drastically when considering complex reflection
groups. For instance, if $W$ has type $G_{13}$ and we choose $c \in \RR/W$
of cardinality 6, then the matrix of $(\ | \ )_c$ for $m=1$ has the
following form.
$$
\begin{pmatrix}
0 & 2 &2 &2 &2 &8 \\ 2 &0 &2 &8 &2 &2 \\ 2 &2 &0 &2 &8 &2 \\ 2 &8 &2 &0 &2 &2 \\ 2 &2 &8 &2 &0 &2 \\
 8 &2 &2 &2 &2 &0
\end{pmatrix}
$$
 
\begin{sidewaystable}
\caption{Discriminants in types $G(e,e,r)$ and exceptional types}
\label{table1}
$$
\begin{array}{|l|c|c|c|}
\hline
e/r & 3 & 4 & 5 \\
\hline
\hline
3 & -(m-9)m^8  & (m-15)(m-3)^3 m^{12}(m+3)^2 & (m-21)(m-6)^4m^{20}(m+3)^5 \\
4 & (m-11)(m-3)^3(m+1)^8  & (m-19)(m-3)^9(m+1)^{12}(m+5)^2 & (m-27)(m-7)^4(m-3)^{10}(m+1)^{20}(m+5)^5 \\
5 & -(m-15)m^{14} & (m-25)(m-5)^3m^{24}(m+5)^2 & (m-35)(m-10)^4m^{40}(m+5)^5 \\
6 & (m-17)(m-5)^3(m+1)^{14} & (m-29)(m-5)^9(m+1)^{24}(m+7)^2 & (m-41)(m-11)^4(m-5)^{10}(m+1)^{40}(m+7)^5 \\
7 & (m-21)m^{20} & (m-35)(m-7)^3 m^{36} (m+7)^2 & (m-49)(m-14)^4m^{60}(m+7)^5 \\
8 & (m-23)(m-7)^3(m+1)^{20} & (m-39)(m-7)^9(m+1)^{36}(m+9)^2 & \\
9 & -(m-27)m^{26} & (m-45)(m-9)^3 m^{48} (m+9)^2 & \\
10 & (m-29)(m-9)^3(m+1)^{26} & & \\
11 & -(m-33)m^{32} & & \\
12 & (m-35)(m-11)^3(m+1)^{32} & & \\
13 & -(m-39)m^{38} & & \\
14 & (m-41)(m-13)^3(m+1)^{38} & & \\
\hline
\end{array}
$$
{}
$$
\begin{array}{l|l|}
 & \mbox{Discriminants} \\
\hline
G_{12} & (m-11)(m-3)^6(m+1)^2(m+5)^3 \\
G_{13} & (m-17)(m-5)^2(m+7)^3 , \ \ \ \  (m-17)(m-5)^2(m-1)^6(m+7)^3\\
G_{22} & (m-29)(m-5)^{15}(m+1)^8(m+11)^6 \\
G_{23} & (m-13)(m-1)^{10}(m+2)^4 \\
G_{24} & (m-17)(m-3)^{12}(m+4)^8 \\
G_{27} & (m-41)(m-5)^{10}(m-1)^{18}(m+4)^{16} \\
G_{28} & (m-15)(m-3)^2(m+1)^9 , \ \ \ \  (m-15)(m-3)^2(m+1)^9 \\
G_{29} & (m-31)(m-11)^4(m+1)^{35} \\
G_{30} & (m-45)(m-5)^{18}m^{16}(m+3)^{25} \\
G_{31} & (m-45)(m-21)^5(m-5)^9(m+3)^{45} \\
G_{33} & (m-33)(m-3)^{24} (m+3)^{20} \\
G_{34} & (m-81)(m-9)^{35}(m+3)^{90} \\
G_{35} & (m-21)(m-3)^{20}(m+3)^{15} \\
G_{36} & (m-33)(m-5)^{27}(m+3)^{35} \\
G_{37} & (m-57)(m-9)^{35}(m+3)^{84} \\
\hline
\end{array}
$$
\end{sidewaystable}

\begin{sidewaystable}
\caption{Discriminants in type $G(2e,e,r)$}
\label{table2}
$$
\begin{array}{|l||c|c|}
e & \multicolumn{2}{|c|}{\mbox{Discriminants for $G(2e,e,3)$}}   \\
\hline
\hline
2 & -(m-9)(m+3)^2  & (m-13)(m-5)^3(m-1)^2(m+3)^6  \\
3 & -(m-13)(m+5)^2 & (m-19)(m-3)^3(m-1)^8(m+3)^6  \\
4 &  -(m-17)(m+7)^2 & (m-25)(m-9)^3(m-1)^8(m+3)^{12}  \\
5 &  -(m-21)(m+9)^2  & (m-31)(m-7)^3(m-1)^{14}(m+3)^{12}  \\
6 & -(m-25)(m+11)^2 & (m-37)(m-13)^3(m-1)^{14}(m+3)^{18} \\
7 & -(m-29)(m+13)^2 & (m-43)(m-11)^3(m-1)^{20}(m+3)^{18} \\
8 & -(m-33)(m+15)^2 & (m-49)(m-17)^3(m-1)^{20}(m+3)^{24} \\
9 & -(m-37)(m+17)^2 & (m-55)(m-15)^3(m-1)^{26}(m+3)^{24} \\
10 & -(m-41)(m+19)^2 & (m-61)(m-21)^3(m-1)^{26}(m+3)^{30} \\
\hline
\end{array}
$$
{}
$$
\begin{array}{|l||c|c|}
e & \multicolumn{2}{|c|}{\mbox{Discriminants for $G(2e,e,4)$}}  \\
\hline
\hline
2 &  (m-13)(m+3)^3 & (m-21)(m-5)^9(m+3)^{14} \\
3 &  (m-19)(m+5)^3 & (m-31)(m-7)^3(m-3)^6(m-1)^{12}(m+3)^{12}(m+5)^2 \\
4 &   (m-25)(m+7)^3 & (m-41)(m-9)^9(m-1)^{12}(m+3)^{24}(m+7)^2 \\
5 &   (m-31)(m+9)^3 & (m-51) (m-11)^3(m-7)^6(m-1)^{24}(m+3)^{24}(m+9)^2 \\
6 & (m-37)(m+11)^3 & (m-61)(m-13)^9(m-1)^{24}(m+3)^{36}(m+11)^2 \\
7 & (m-43)(m+13)^3 & (m-71)(m-15)^3(m-11)^6(m-1)^{36}(m+3)^{36}(m+13)^2 \\
8 & (m-49)(m+15)^3 & (m-81)(m-17)^9(m-1)^{36}(m+3)^{48}(m+15)^2 \\
\hline
\end{array}
$$
{}
$$
\begin{array}{l||c|c|c|}
e & \multicolumn{3}{|c|}{\mbox{Discriminants for $G(2e,e,2)$}} \\
\hline
1 & (m-3)(m+1) & (m-3)(m+1) & \\
2 & (m-5)(m+3) & (m-5)(m+3) & (m-5)(m+3) \\
3 & (m-7)(m+5) & (m-7)(m-3)(m-1)^2(m+3)^2 & \\ 
4 & (m-9)(m+7) & (m-9)(m-1)(m+3)^2 & (m-9)(m-1)(m+3)^2 \\
5 & (m-11)(m+9) & (m-11)(m-7)(m-1)^4(m+3)^4 & \\
6 & (m-13)(m+11) & (m-13)(m-1)^2(m+3)^3 & (m-13)(m-1)^2(m+3)^3 \\
7 & (m-15)(m+13) & (m-15)(m-11)(m-1)^6(m+3)^6 & \\
8 & (m-17)(m+15) & (m-17)(m-1)^3(m+3)^4 & (m-17)(m-1)^3(m+3)^4 \\
9 & (m-19)(m+17) & (m-19)(m-15)(m-1)^8(m+3)^8 & \\
10 & (m-21)(m+19) & (m-21)(m-1)^4(m+3)^5 & (m-21)(m-1)^4(m+3)^5 \\
\hline
\end{array}
$$
\end{sidewaystable}

\end{document}